\documentclass{imsart}

\RequirePackage[OT1]{fontenc}

\usepackage{xr-hyper}
\usepackage{hyperref}

\usepackage{amsmath,amsthm,amssymb}
\usepackage[utf8]{inputenc}
\usepackage{amsfonts,dsfont}
\usepackage{soul,graphicx}
\graphicspath{
  {fig/}
}
\usepackage{natbib}

\startlocaldefs

\usepackage{amsthm}
\newtheorem{theorem}{Theorem}[section]
\newtheorem{proposition}[theorem]{Proposition}

\newtheorem{lemma}[theorem]{Lemma}

\theoremstyle{remark}
\newtheorem{remark}[theorem]{Remark}

\theoremstyle{definition}
\newtheorem{definition}[theorem]{Definition}

\usepackage[ruled]{algorithm2e}

\newcommand{\R}{\mathbb{R}}
\newcommand{\Rfam}{\mathfrak{R}}

\renewcommand{\P}{\mathbb{P}}
\newcommand{\E}{\mathbb{E}}

\newcommand{\telque}{\,:\,}

\newcommand{\kFWER}{\mbox{$k$-FWER}}

\newcommand{\Pow}{\mbox{Pow}}

\newcommand{\JER}{\mathrm{JER}}
\newcommand{\JC}{\mathcal{E}}

\newcommand{\mbn}{\mathbb{N}}

\newcommand{\mtc}{\mathcal}

\newcommand{\wt}[1]{{\widetilde{#1}}}
\newcommand{\wh}[1]{{\widehat{#1}}}
\newcommand{\ol}[1]{\overline{#1}}
\newcommand{\ind}[1]{{\mathds{1}\left\{#1\right\}}}

\newcommand{\paren}[1]{\left(#1\right)}

\newcommand{\set}[1]{\left\{#1\right\}}
\newcommand{\abs}[1]{\left| #1 \right|}

\newcommand{\cA}{{\mtc{A}}}

\newcommand{\cH}{{\mtc{H}}}

\newcommand{\cX}{{\mtc{X}}}
\newcommand{\cP}{{\mtc{P}}}

\newcommand{\cN}{{\mtc{N}}}

\newcommand{\magen}[1]{
#1
}

\newcounter{nbdrafts}
\makeatletter
\newcommand{\checknbdrafts}{
\ifnum \thenbdrafts > 0
\@latex@warning@no@line{**********************************************************************}
\@latex@warning@no@line{* The document contains \thenbdrafts \space draft note(s)}
\@latex@warning@no@line{**********************************************************************}
\fi}

\makeatother

\newcommand{\Nm}{\mathbb{N}_m}
\newcommand{\Nn}{\mathbb{N}_n}

\newcommand{\iop}{{post hoc }}
\newcommand{\iopeq}{$\mathrm{PH}_\alpha$}

\endlocaldefs

\begin{document}

\begin{frontmatter}

\title{Post hoc inference via joint family-wise error rate  control}
\runtitle{Post hoc inference via JER control}
\begin{aug}
\author{\fnms{Gilles} \snm{Blanchard}
\ead[label=e1]{gilles.blanchard@math.uni-potsdam.de}}
\address{
Universit\"at Potsdam, Institut f\"ur Mathematik\\
Karl-Liebknecht-Stra{\ss}e 24-25 14476 Potsdam, Germany\\
\printead{e1}\\
\phantom{E-mail: gilles.blanchard@math.uni-potsdam.de\ }}
\and
\author{\fnms{Pierre} \snm{Neuvial}\ead[label=e2]{pierre.neuvial@math.univ-toulouse.fr}}
\address{
Institut de Mathématiques de Toulouse; \\
UMR 5219, Université de Toulouse, CNRS\\
UPS IMT, F-31062 Toulouse Cedex 9, France\\
\printead{e2}\\
\phantom{Email: pierre.neuvial@math.univ-toulouse.fr\ }}
\and
\author{\fnms{Etienne} \snm{Roquain}
\ead[label=e3]{etienne.roquain@upmc.fr}}
\address{
Sorbonne Universit\'e (Universit\'e Pierre et Marie Curie), LPSM,\\ 4, Place Jussieu, 75252 Paris cedex 05, France\\
\printead{e3}\\
\phantom{E-mail: etienne.roquain@upmc.fr }}

\end{aug}
\begin{abstract}
  We introduce a general methodology for post hoc inference in a large-scale multiple testing framework.  The approach is called ``user-agnostic'' in the sense that the statistical guarantee on the number of correct rejections holds for any  set of candidate items selected by the user (after having seen the data).  This task is investigated by defining a suitable criterion, named the joint-family-wise-error rate (JER for short). We propose several procedures for controlling the JER, with a special focus on incorporating  dependencies while adapting to the unknown quantity of signal (via a step-down approach).  We show that our proposed setting incorporates as particular cases a version of the higher criticism as well as the closed testing based approach of \citet{GS2011}. Our theoretical statements are supported by numerical experiments.
\end{abstract}



\end{frontmatter}


\section{Introduction}
\label{sec:context}

\begin{figure}[htp!]
  \centering
 \includegraphics[scale=0.3]{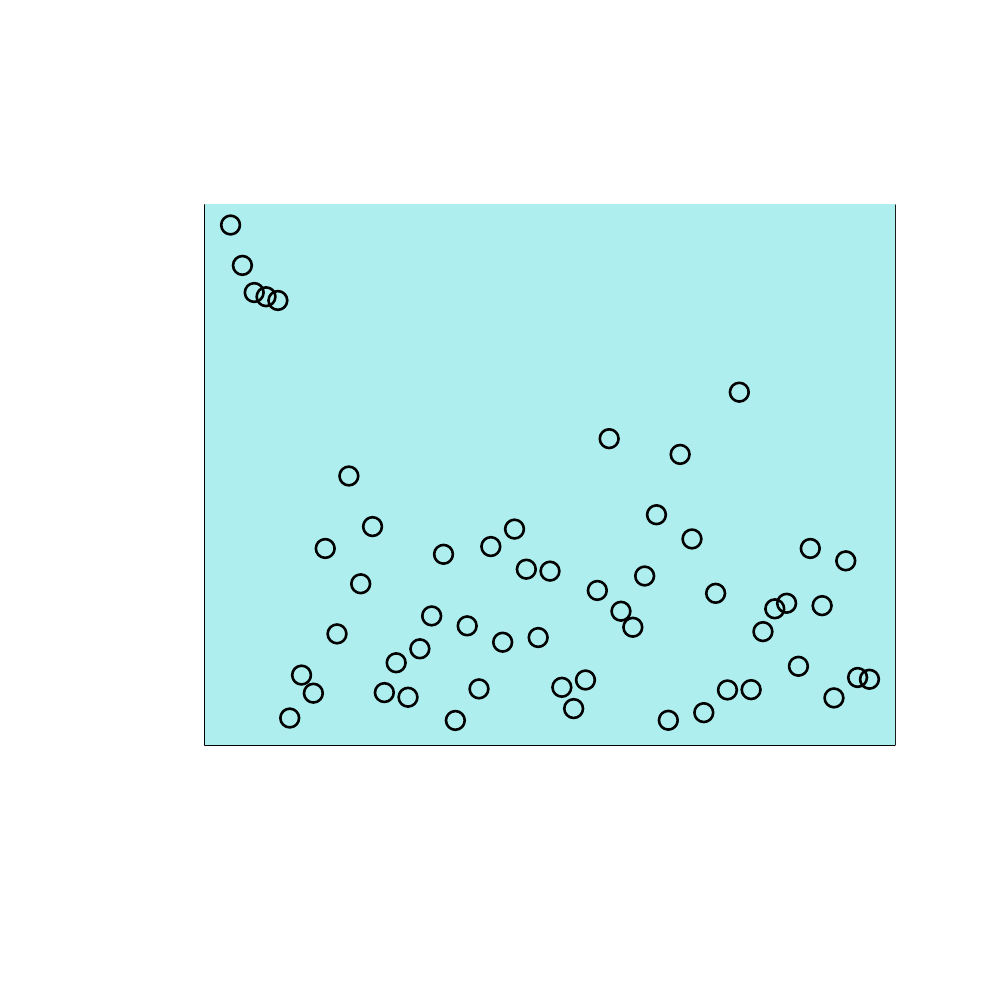} \includegraphics[scale=0.3]{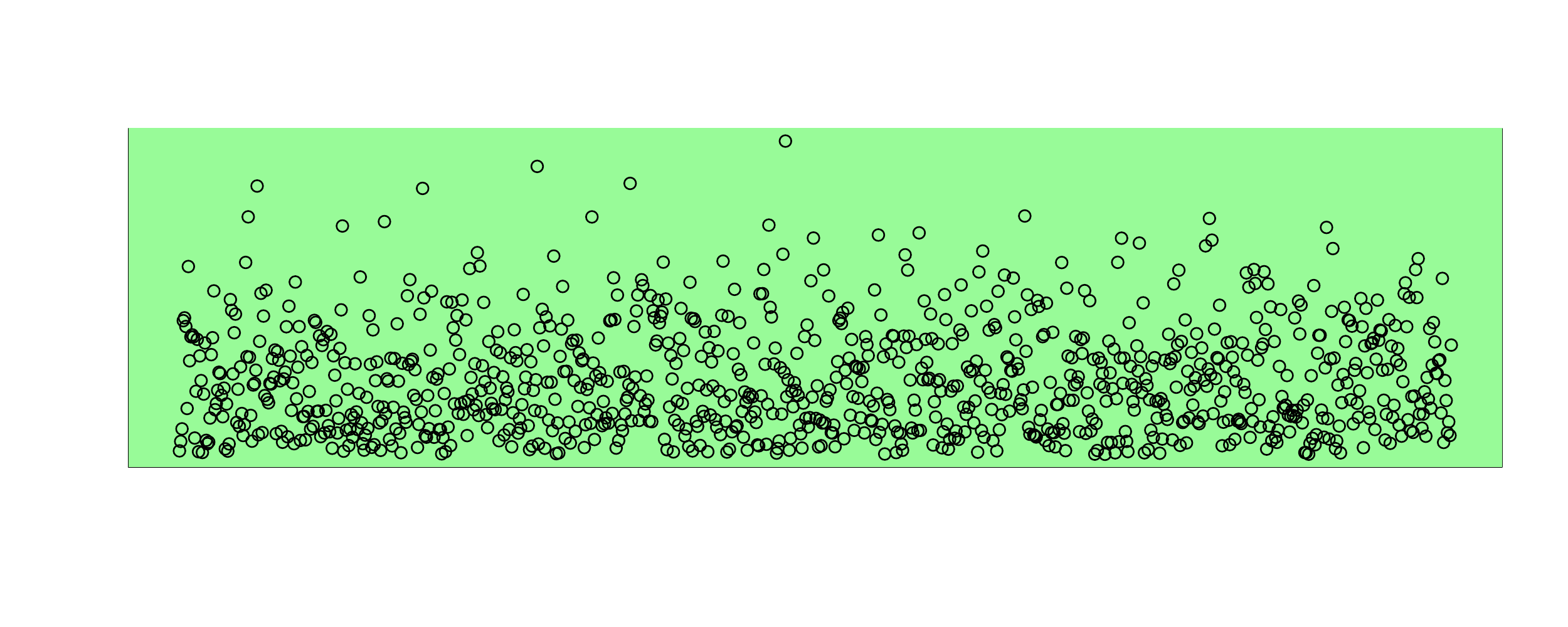}
\caption{Illustration of the post hoc selection effect.  Right: virtual data set with $1000$ measurements.
Left: data set of $55$ measurements  selected  from the right dataset.
Measures have been generated as i.i.d. absolute values of $\mathcal{N}(0,1)$. }  \label{fig:selectioneffect}
\end{figure}

Large-scale multiple inference with a rigorous statistical guarantee
has become a topic of ever increasing relevance with the advent of very high-dimensional
data in numerous application areas.
Classical multiple testing procedures prescribe a rejection set based on the amount of false positives that the user might tolerate (e.g., false discovery rate control at level $5\%$). However, if the result does not correspond to what the user expected,
they may tend to ``snoop'' in the data, possibly concentrating only on a set $R$ of hypotheses that appear promising to them.
Even when motivated by plausible justifications, any such approach
will invalidate standard statistical guarantee because of the {\it selection effect}.
This is illustrated on Figure~\ref{fig:selectioneffect}, where only ``noisy" measurements have been generated:
within the selected set (in blue), $5$ points look like significant measurements. However, this is only due to the selection effect: the blue data set comes from a larger data set (green) where these $5$ measures are just the $5$ maximum (noisy) measurements.
As a consequence, while building a statistical guarantee on the selected set $R$,  the overall size of the data set should be considered.
This is the aim of the so-called ``post-selection'' (or post hoc) inference.

A particular case of post hoc inference is faced when the selection step $R$ is a pre-specified selection method,  see \cite{BY2005,Loc2013,FST2014,BM2014,BCH2014,TT2015}, among others. 
However, since the selection step is fixed, this does not allow for arbitrary ``data snooping'' or {\em ad hoc} selection rules often used in exploratory research.

More generally, elaborate selection rules possibly consisting in several stages and
involving user-fixed tuning constants are commonly used in a variety of contexts, for instance:
\begin{itemize}
\item in neural activity detection from brain imaging data, cluster-extent approaches \citep{woo14cluster} select voxels by a two-stage process, first building groups of contiguous voxels whose activity levels all pass a user-defined threshold, then performing a correction to select a subset of clusters. The second stage only ensures that each cluster contains at least one truly active voxel, but there is no additional statistical
guarantee about the proportion of active voxels among the selected.
\item in the context of gene or protein activity change detection,
  a two-sample rank test might
  be used to detect activity changes, while requiring that the log-ratio of average
  observed activities of the two samples (``fold change'') is larger than a certain user-specified level, see \cite{li2012volcano}. 
In other words, for each hypothesis a statistic $T_1$ is used for constructing a standard test, but a different statistic $T_2$ is used for screening, with the two statistics not being independent. 
\end{itemize}

A point of view argued in several papers in various statistical contexts, e.g., \cite{GS2011,Ber2013,BPS2016} is that in absence of precise information of the
user's selection strategy, it is desirable to  provide a statistical guarantee {\it simultaneously} for any possible selected set.
In this paper, we adopt this view and focus on  simultaneous upper bounds on the number of false positives on the selected set, as proposed in the seminal paper \cite{GS2011}.
More formally, our goal is to build functional $V(\cdot)$ defined on all subset of hypotheses, such that the following uniform guarantee holds:
\label{eq:easy-JER}
$$
\P( \forall R\subset \{1,\dots,m\}\::\: |\cH_0\cap R|\leq V(R))\geq 1-\alpha,
$$
where $m$ is the number of null hypotheses to be tested (identified with their respective index) and $\cH_0\subset \{1,\dots,m\}$ corresponds to the (unknown) set of true null hypotheses.
 This general principle is ``user-agnostic'', in the sense that the provided inference is ``ready for any selected set'' (the ``for all $R$'' being inside the probability).
 Observe that a bound $V(\cdot)$ satisfying the above guarantee can also inform the
 choice of the final rejected set $R$; for example the user is allowed to optimize
 some function of $V(R)$, possibly subject to geometrical or data-dependent constraints on $R$.

Our construction of post hoc bounds relies on the control of a multiple testing criterion that we call ``joint (family-wise) error rate'' (JER for short), which was implicitly defined by \citet{Mein2006} for building false discovery proportion confidence envelopes (see also \citealp{GW2004,GW2006} for more details on this topic).
The JER has a particularly simple expression in the case of $p$-value thresholding:
given a family $\{p_i(X),1\leq i\leq m\}$ of $m$ $p$-values and a family of thresholds $\mathcal{T}=(t_k)_{1\leq k\leq K}$, the JER of $\mathcal{T}$ is related to the distribution of $p_{(k:\cH_0)}$, the $k$-th smallest value in the set $\{p_i,i\in \cH_0\}$:
\begin{equation}\label{equ-JR4pvalue}
\JER(\mathcal{T})=
\P\bigg(
\exists k \in\{ 1,\dots,K\wedge m_0\}\: :\:
p_{(k:\cH_0)} < t_k
 \bigg),
\end{equation}
where $m_0=|\cH_0|$ is the number of true null hypotheses.
It turns out that finding $\mathcal{T}$ such that $\JER(\mathcal{T})\leq \alpha$ provides that the functional
\begin{align}
{V}(R) &= \min_{k \in \{1, \dots, K\}}\left \{\sum_{ i \in R} \ind{p_i(X)\geq  t_k} + k-1\right\} 
\:,\:\: R\subset \{1,\dots,m\}\:\label{posthocV}
\end{align} 
is a valid post hoc bound  (see Section~\ref{se:JRglobal} for a proof in a  general context).
Hence, a general intuition is that the threshold $t_k$ should be chosen as an appropriate quantile of the distribution of $p_{(k:\cH_0)}$, with some extra slack to take into account for uniformity in $k$.
\\

The contributions of the present work are the following:
\begin{itemize}
\item a general framework to build post hoc bounds, that generalizes the method of \cite{GS2011} and does not rely on closed testing but on JER control;
\item JER controlling procedures, with adaptivity to dependence and to the proportion of true null hypotheses. 
  These procedures are implemented in an open-source \textsf{R}~\citep{R} package \citep{sanssouci-package};
\item reproducible numerical experiments to illustrate our theoretical statements.
\end{itemize}

In addition, this study connects former ({\it a priori} unrelated) concepts:
the closed testing-based method of \cite{GS2010,GS2011},
the confidence envelopes of \cite{Mein2006} and the higher criticism of \cite{DJ2004}.

The paper is organized as follows. In Section~\ref{se:JRglobal}, we expose the general approach to post hoc
multiple test inference based on JER control. In the following sections, we develop this point of view in some specific exemplary models under known or unknown dependence structure; the models are presented in Section~\ref{sec:modelassumption} and the basic JER control obtained using the classical Simes inequality
is analyzed in Section~\ref{sec:basicineq}.
In Section~\ref{sec:improving-jr-control}, we present
improvements to this basic case by considering more general threshold families,
incorporating adaptation to noise dependence structure, and
a step-down principle. Two specific examples of this
improved methodology are developed in Section~\ref{sec:two-templates}. In Section~\ref{sec:num}, we present the results of numerical simulations
illustrating and comparing the developed methods. We conclude with a discussion of
various points in Section~\ref{se:literature}.
Due to space constraints, proofs as well as some additional results (including a detailed
comparison to the work of \citealp{GS2010,GS2011}, to the higher criticism of \citealp{DJ2004},
related optimality properties for detection purposes, and algorithmic details concerning
Monte-Carlo and permutation-based calibration)
are postponed to the supplementary material \cite{BNR2017supp}.
The sections of this supplement are referred to with an additional symbol ``S-'' in the  numbering.

\section{JER control: principle and properties}
\label{se:JRglobal}

In this section, we introduce the framework (Section \ref{sec:aim}) for post hoc multiple testing inference,
and propose a general approach to tackle this problem based on a reference family of rejection sets (Section \ref{sec:jfweintro}). Proceeding from the general to the particular, we will first study and discuss some generic properties of this approach (Section \ref{sec:propnewbounds}) before focusing
on more specific choices for the reference family leading to \eqref{equ-JR4pvalue} and \eqref{posthocV} (Section \ref{se:thresholdbased}). 
Formal proofs for theoretical claims in this section are found in Section~\ref{app:proofsgeneric}.
\subsection{Aim}\label{sec:aim} 

Formally, let $X$ denote observed data generated from a statistical
model $(\cX,\mathfrak{X},P)$, $P\in \mtc{P}$, and assume we want to
test a collection of null hypotheses $H_{0,i} \subset \mtc{P}$
indexed by $i \in \Nm := \set{1,\ldots,m}$.  For any $P\in\mtc{P}$, we
denote by $\cH_0(P)$ the set of (indices of) true null hypotheses satisfied by $P$,
that is, $\cH_0(P)=\{i \in \Nm \telque P \in H_{0,i}\}$, and by
$m_0(P)$ its cardinality (or $\cH_0$, $m_0$ for short). We denote by $\pi_0=m_0/m$ the proportion of true nulls. We also let $\cH_1(P)=\Nm \backslash \cH_0(P)$ be the set of (indices of) false nulls and $m_1(P)=m-m_0(P)$ its cardinality (or $\cH_1$, $m_1$ for short).

Our main objective in this paper is to find a function ${V}(X,R)$ (denoted by ${V}(R)$ for short) satisfying
\begin{equation}\label{aim}\tag{\iopeq}
\text{for all } P\in\mtc{P},\, \qquad \P_{X\sim P}\bigg(\forall R \subset 
\Nm ,\:
|R\cap \cH_0(P)| \leq {V}(R) \bigg)\geq 1-\alpha.
\end{equation}
If the above is satisfied, $V(R)$ gives  a level $1-\alpha$ confidence bound for the number of false rejections
in a set $R$ of (indices of) rejected hypotheses that is {\em uniformly valid} over all possible choices of $R$.
Letting ${S}(R)= |R| - {V}(R)$, the property \eqref{aim} equivalently provides the following simultaneous lower bound on $|R\cap \cH_1(P)|$, that is, evidence of signal in $R$:
\begin{equation*}
\text{for all } P\in\mtc{P},\, \qquad \P_{X\sim P}\bigg(\forall R \subset 
\Nm ,\:
|R\cap \cH_1(P)| \geq {S}(R) \bigg)\geq 1-\alpha\,.
\end{equation*}
As the the above bounds are uniformly valid over all possible choice of $R$, they will apply (with probability at least $1-\alpha$) to any arbitrary data-dependent choice of $R$
made by the user, including choices made after looking at the value of the bound itself for different candidates for $R$. For instance, $R$ can be chosen as maximizing $|\hat{R}|$ among those $\hat{R}$ satisfying $S(\hat{R})/|\hat{R}|\geq 0.5$ (more than half of signal in $\hat{R}$ with high probability). Obviously, the theoretical guarantees for $\hat{R}$ also hold because the bounds are uniform on $R\subset \Nm$.

\subsection{General principle}
\label{sec:jfweintro}

The question of how to obtain a control of the general form \eqref{aim} is statistical
as well as computational in nature, since it is not practically feasible to consider individually all $2^m$ possibilities for candidate rejection sets $R$ as soon as $m$ exceeds a couple of dozens. Provided that the statistical guarantee holds, we would ideally wish that the bound ${V}(R)$ is computable efficiently for any candidate $R$ (or family thereof) suggested by the user.

\begin{figure}[!htp]
  \centering
    \includegraphics[scale=0.3]{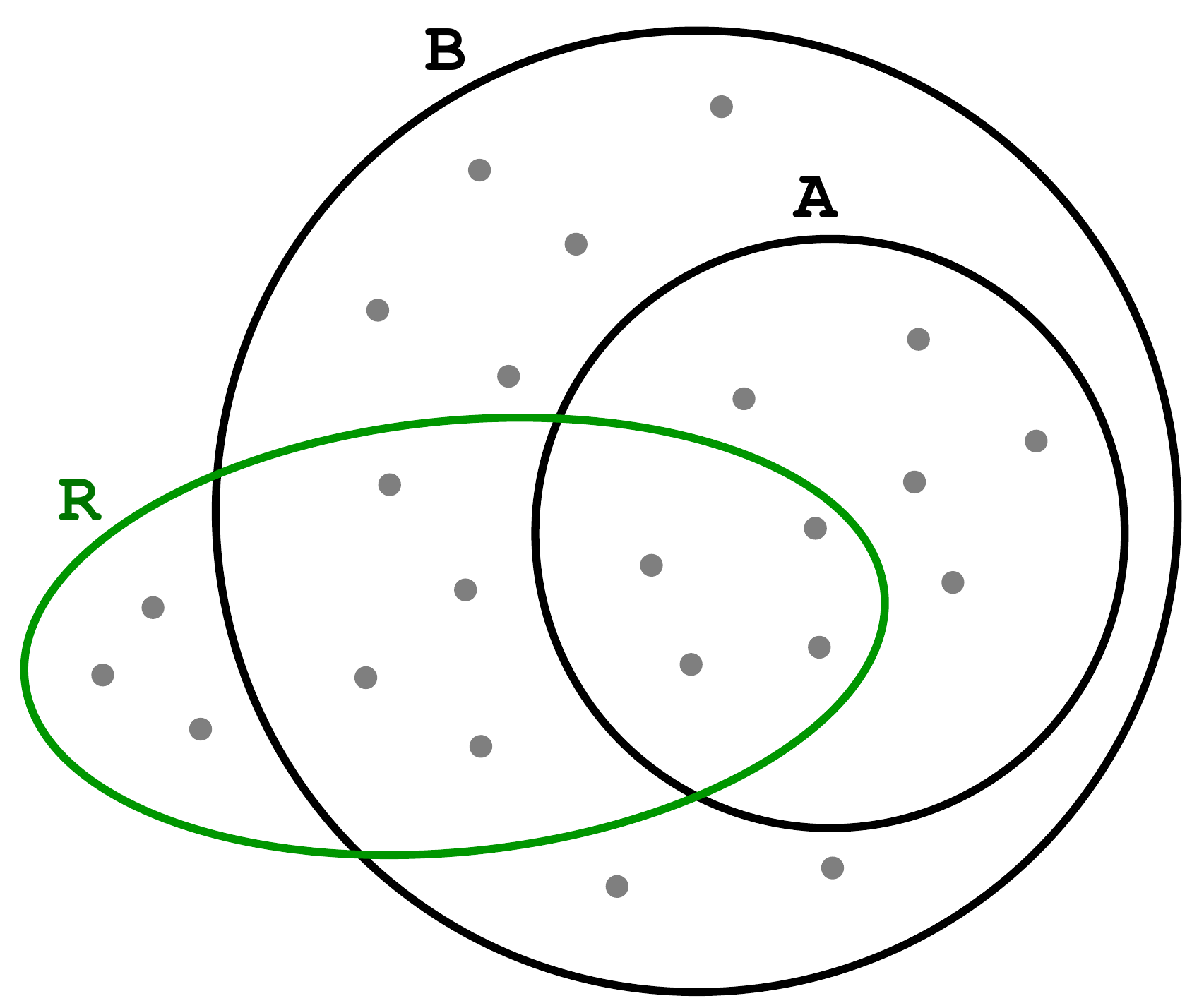}
  \caption{Toy example: use of a reference family with two subsets $A$ and $B$ to build a post hoc bound on the number of true positives in an arbitrary candidate rejection set $R$.}
  \label{fig:illustration-bound}
\end{figure}

In this section, we consider a general approach to the problem based on a reference family with a controlled Joint family-wise Error Rate (JER). The basic argument is illustrated by Figure~\ref{fig:illustration-bound}. Imagine that a subset $A$ of hypotheses is guaranteed to contain less than 5 true nulls, that is, $|A\cap \cH_0(P)|\leq 5$. Then this also provides information on other subsets $R\subset \Nm$ with $R\neq  A$. Namely, for any $R\subset \Nm$, $|R\cap \cH_1(P)|\geq |R\cap A|-5$. Of course, while this information is useful for $R$ if $|R\cap A|\geq 6$, it is not if $|R\cap A|\leq 5$ (nonpositive bound), as in Figure~\ref{fig:illustration-bound}. Next, if we want  to improve the bound, we can consider another set $B$ (here including $A$) with the property
$|B\cap \cH_0(P)|\leq 7$ (say). In the situation pictured in Figure~\ref{fig:illustration-bound}, this ensures that $R$ contains at least one element which is in $\cH_1(P)$.

More generally, let us assume that we have at hand $\Rfam=((R_1(X),\zeta_1(X)),$ $\ldots,$ $(R_K(X),\zeta_K(X)))$ a data-dependent collection of
subsets $R_k$ of $\Nm$ and integer numbers $\zeta_k$ (we will often omit the dependence in $X$ to ease notation), such that, with probability larger than $1-\alpha$, the set $R_k(X)$ does not contain more than $\zeta_k(X)$ elements of $\cH_0(P)$, uniformly over $k$, that is,

\begin{equation}
\text{For all } P \in \cP, \qquad
\JER(\Rfam,P)\leq \alpha,
\label{eq:jfwer}
\end{equation}
where we have denoted
\begin{equation}
\JER(\Rfam,P) := 1-
\P_{X\sim P}(\JC(\Rfam,\cH_0(P)))\,,
\label{eq:jr}
\end{equation}
with the event
\begin{equation}
 \JC(\Rfam,\cH_0) := \set{\forall k = 1,\dots,K,\;
|R_k(X)\cap \cH_0| \leq \zeta_k(X)}.
\label{eq:jrdef}
\end{equation}

We see $\Rfam$ as a {\em reference family} of rejection sets for which a statistical guarantee on the number of false rejections is ensured, and based on which we will build a \iop bound. The cardinality (or size) $K$ of the reference family is also allowed to be data-dependent in the most general form, although this dependence is not acknowledged for in our notation for simplicity. Different choices are possible for $\Rfam$, allowing
to recover as particular cases settings considered in previous literature.
Let us mention two important cases concerning the bounds $\zeta_k$:
\label{page-seven}
\begin{itemize}
\item $\zeta_k = k-1$ for all $k$: in this case, each individual rejection region $R_k$ has controlled $k$-FWER, and the control is uniform over the regions.
\item $\zeta_k = |R_k|-1$ for all $k$:  adopting a different point of
view, associate to each $R\subset \Nm$
the {\em intersection hypothesis} $H_{0,R}:= \bigcap_{i \in R} H_{0,i}$ (in this view, each $R$ corresponds to a hypothesis rather than a collection of hypotheses).
The statement \eqref{eq:jfwer}-\eqref{eq:jrdef} is interpreted
as saying that with high probability, each individual rejection region $R_k$ has at least one true rejection. Consequently, rejecting all intersection hypotheses $H_{0,R_k}$, $k=1\dots,K$ can be done without committing any error. This corresponds to an overall family-wise error rate control over this family of hypotheses. 
\end{itemize}

Our first goal in this paper is to analyze how to go from the JER control~\eqref{eq:jfwer}-\eqref{eq:jrdef} to a post hoc statement \eqref{aim}. This will be done in the present section in a general setting. In the remaining sections,
 we will concentrate on how to obtain the JER control itself.
 For this, we will focus on the first situation above ($\zeta_k=k-1$)
 and therefore assume this setting by default unless otherwise specified.
 In the second situation ($\zeta_k=|R_k|-1$),
 JER control
 can in particular be obtained via closed testing, thus recovering the setting
 of \cite{GS2011}, see Section~\ref{sec:compclosed} for a more detailed discussion.

How can we ``interpolate'' from the control on a reference family \eqref{eq:jfwer} to a control on all possible rejection sets \eqref{aim}?
On the event \eqref{eq:jrdef}, the only available information on the unknown subset $\cH_0$ is that it is an element of the collection of subsets
\begin{align*}
\cA(\Rfam) &:= \set{A\subset \Nm: \JC(\Rfam,A) \text{ holds }} \\
           &\:= \set{A \subset \Nm:\;\forall k=1,\dots,K,\;
|R_k\cap A| \leq \zeta_k}.
\end{align*}
As a result, the best we can do to bound $|R \cap \cH_0|$
for any proposed rejection set $R$ is a worst-case bound under this constraint:
\begin{equation}\label{eq:optbound1}
  V^*_\Rfam(R) := \max_{A \in \cA(\Rfam)} \abs{R\cap A}, \:\:\: R\subset \Nm\,.
\end{equation}
A significant problem is that $V^*(R)$ (we will sometimes drop the index $\Rfam$ for simplicity)
may not be easy to compute in general (see Proposition~\ref{prop:NPhard} below). 
We therefore introduce the following coarser but simpler bound:
\begin{equation}
\ol{V}_\Rfam(R) :=
\min_{k \in \{1, \dots, K\}} \paren{ \abs{R\setminus R_k} + \zeta_k} \wedge |R|\,,\:\: R\subset \Nm\,.
\label{eq:aug}
\end{equation}
Observe that $\ol{V}(R)$ is non-decreasing in the sense that
$R\subset R'$ implies $\ol{V}(R)\leq \ol{V}(R')$. The next result formalizes the link between JER control and the associated \iop bounds.
This result is proved in Section~\ref{proof:se:JRglobal}, along with all of the other results of the present section.

\begin{proposition}\label{prop:Vbarisposthoc}
  Let $\Rfam=(R_k(X),\zeta_k(X))_{1\leq k\leq K}$ be a data-dependent collection of subsets $R_k$ of $\Nm$ and of integers $\zeta_k$\,. Then for any $\cH_0 \subset \Nm$, $\cH_1=\Nm \setminus \cH_0$\,, the event  $\JC(\Rfam,\cH_0)$ defined in \eqref{eq:jrdef} is such that
\begin{align}
 \JC(\Rfam,\cH_0)& =
\left\{\forall R \subset \Nm ,\:
|R\cap \cH_0| \leq \ol{V}_\Rfam(R)  \right\}
\label{equ:event} \\
& = \left\{\forall R \subset \Nm ,\:
|R\cap \cH_0| \leq V^*_\Rfam(R)  \right\}
 \label{equ:eventopt} \,.
\end{align}
\end{proposition}

In particular, Proposition~\ref{prop:Vbarisposthoc} shows that $\Rfam$ satisfies the JER control \eqref{eq:jfwer} if and only if
$\ol{V}_\Rfam(\cdot)$ or $V^*_\Rfam(\cdot)$ satisfies \eqref{aim}.

\subsection{General properties}\label{sec:propnewbounds}
In this section, we further discuss general properties of the obtained post hoc bounds. 
The JER control gives rise to the \iop upper-bound $\ol{V}$, which we can see as an approximation of the optimal bound $V^*$\,.
A first legitimate question is whether an approximation of the optimal bound is necessary in the first place,
and then whether these approximations possess favorable properties.
In this section, we provide arguments in this direction.

\begin{remark}
\label{rmk:equiv-V-S}
The results of the paper can equivalently be stated in terms of false positives using $V$, $V^*$ and $\ol{V}$ or in terms of true positives $S$, $S^*$ and $\ol{S}$, where for any $R \in \Nm$ $S^*(R):=|R|-V^*(R)$ and $\ol{S}(R) := |R| - \ol{S}(R)$. For simplicity we have chosen to focus on $V$. 
\end{remark}

\paragraph{Computing the optimal bounds is NP-hard}

The claim that computing the optimal bound $V^*$ is computationally difficult in general is supported by the following NP-hardness result:
\begin{proposition}
  \label{prop:NPhard}
  The problem of computing $V_\Rfam^*(R)$
given an arbitrary reference family $\Rfam=(R_k,\zeta_k)_{1\leq k\leq K}$ (with $R_k\subset \Nm,\zeta_k \in \mbn$\,),
and $R\subset \Nm$\,, is NP-hard.
\end{proposition}

Naturally, Proposition~\ref{prop:NPhard} does not imply that computing
the optimal bound $V^*(R)$ is always infeasible: depending on the
choice of the reference family, we might be in a particular case where
this can be done efficiently --- in fact, we will discuss precisely
such a situation below (nested regions).
Still, it is worth
noting that the proof of the above result establishes NP-hardness for
the more specific case $\zeta_k = |R_k|-1$\,, where the reference
family is interpreted as tests of certain intersection hypotheses.  We
show in Section~\ref{sec:compclosed} that in this case, the bound
$V^*$ coincides with the bound that can be derived from the closed testing approach
of \citet{GS2011}.

In general, it is therefore sensible in practice to look for computable approximations of $V^*$. We turn to general properties of the proposed bound $\ol{V}$\,.

\paragraph{Self-consistency}

Given some reference family $\Rfam=(R_k,\zeta_k)_{1\leq k \leq K}$, on the large probability event \eqref{equ:event} for which the control $|R_k\cap \cH_0(P)| \leq \zeta_k$, $1\leq k \leq K$ holds, $\ol{V}_\Rfam$ provides a bound for $|R_k\cap \cH_0(P)|$ itself, namely
\begin{equation}\label{equ:zetatilde}
\tilde{\zeta}_k :=  \ol{V}_\Rfam(R_k) = \min_{j \in \{1, \dots, K\}} \paren{ \abs{R_k\setminus R_j} + \zeta_j} \wedge |R_k| ,\:\:\:\:1\leq k \leq K.
\end{equation}
Obviously, $\tilde{\zeta}_k\leq \zeta_k$, with a possible strict inequality.
Nevertheless, the next proposition shows that there is no advantage in ``iterating'' the post hoc bound $\ol{V}$ with $\zeta$ replaced by $\tilde{\zeta}$.
\begin{proposition}\label{prop:iterateisnotmagic}
For any reference family $\Rfam=(R_k,\zeta_k)_{1\leq k \leq K}$, define $(\tilde{\zeta}_k)_{1\leq k \leq K}$ by \eqref{equ:zetatilde}. Denoting $\wt{\Rfam} = (R_k,\wt{\zeta}_k)_{1\leq k \leq K}$, we have 
\begin{align}\label{equ:iterateisnotmagic}
\ol{V}_\Rfam(R) &=
\min_{k \in \{1, \dots, K\}} \paren{ \abs{R\setminus R_k} + \tilde{\zeta}_k} \wedge |R|= \ol{V}_ {\wt{\Rfam}}(R)\,,\:\: R\subset \Nm\, .
\end{align}
\end{proposition}

In particular, the $\tilde{\zeta}_k$s satisfy the following ``self-consistency'' equation:
\begin{equation}\label{equ:equzetatilde}
\tilde{\zeta}_k =  \min_{j \in \{1, \dots, K\}} \paren{ \abs{R_k\setminus R_j} + \tilde{\zeta}_j} \wedge |R_k| ,\:\:\:\:1\leq k \leq K.
\end{equation}

\paragraph{Optimality under nestedness assumption}
In the situation where the sets $(R_k)_{1\leq k \leq K}$ are nested, it holds that $\ol{V}=V^*$, that is, the formula for $\ol{V}$ provides a computationally efficient way to compute the optimal bound in this case.

\begin{proposition}\label{prop:optimalbound}
  For any reference family $\Rfam=(R_k,\zeta_k)_{1\leq k\leq K}$
  such that $R_k\subset R_{k'}$ whenever $k\leq k'$, we have $\ol{V}_\Rfam(R) = V^*_\Rfam(R)\,.$
\end{proposition}

The more specific reference families studied in the remainder of the paper will
satisfy the nestedness assumption, but it is in general not the case for closed testing-based families.

\subsection{Focus of the paper}\label{se:thresholdbased}

The aim of the rest of the paper is to find suitable reference families $\Rfam$ (which may be seen as ``procedures'') that control the joint family-wise error rate at some pre-specified level $\alpha$.

A variety of choices are possible for the reference family. In this paper, we focus on the common situation where a test statistic $T_i(X)$ is available for each null hypothesis $H_{0,i}$,  which in turn is transformed into a $p$-value $p_i(X)$, for all $i \in \Nm$.
As announced earlier, we will also always choose $\zeta_k=k-1$, $1\leq k \leq K$ from now on 
and therefore omit the $\zeta$ and use the simplified notation $\Rfam=(R_1(X),\ldots,R_K(X))$  for the reference family. 
We will also assume that $K$ is  non-random and has been fixed in advance. In this situation, a simple way to build a reference family  is to use $p$-value thresholding:
\begin{equation}\label{equRkthresholdpvalues}
R_k(X)=\{i\in\Nm \::\: p_i(X) < t_k\}\,,\:\:k\in\{1,\dots,K\},
\end{equation}
where the $t_k\in\R$, $1\leq k \leq K$, are associated thresholds, possibly depending on $X$. We easily check that the simpler expressions \eqref{equ-JR4pvalue} and \eqref{posthocV} announced in the introduction hold in that context.

\section{Model assumptions}\label{sec:modelassumption}

Properties of the $p$-value process $(p_i(X), i \in \Nm)$ depend on the underlying model assumptions. In this paper, we distinguish between two general situations, depending on whether the dependence structure is known or not.

\subsection{Location model}\label{sec:leadingmodel}

To give some intuition behind the general assumptions of the next section, we start by considering a specific location model 
\begin{equation}\label{model-highdim}
X_i = \mu_i + \varepsilon_i, \:\:\:i\in\Nm\,,
\end{equation}
where the $\varepsilon_i$ are identically distributed with a common known marginal distribution which is assumed to be continuous, integrable and symmetric.
We denote $\ol{F}(x)=\P(\varepsilon_1 \geq x)$, $x\in \R$.
We consider the one-sided (resp. two-sided) testing problem with null hypotheses $H_{0,i} :$ ``$\mu_i\leq 0$'' (resp. $H_{0,i} :$``$ \mu_i= 0$'') versus the alternative hypotheses $H_{1,i} : $``$\mu_i >0$'' (resp. $H_{1,i} :$``$ \mu_i\neq  0$'') for all $i\in \Nm$.
Classical $p$-values are then given by $p_i(X)=\ol{F}(X_i)$ (resp. $p_i(X)=2\ol{F}(|X_i|)$).
As many procedures of multiple testing theory, our results will rely on the (joint) distribution of
$
(p_i(X))_{i\in\cH_0(P)}
$
or some approximation/bound of it.

\paragraph{Known dependence}
In the case where the (joint) distribution of $\varepsilon$ is known, we can consider ``least favorable" $p$-values $q_i(X)=\ol{F}(X_i-\mu_i)$ ($q_i=2\ol{F}(|X_i-\mu_i|)$. While the $q_i(X)$'s are not observed, they can be used purely as a technical device. Interestingly, these variables satisfy the following point-wise property: for all $i\in\cH_0$, $p_i(X)\geq q_i(X)$, both in the one-sided and two-sided case. In addition, their joint distribution, that is,
$
\nu_m = \mathcal{D}((q_i(X))_{1\leq i \leq m}),
$
is assumed to be known. For instance, under independence of the $\varepsilon_i$'s,
$\nu_m=U(0,1)^{\otimes m}$.
 
\paragraph{Unknown dependence} In the case where the (joint) distribution of $\varepsilon$ is  unknown, so is $\nu_m$ and the above least favorable $p$-values cannot be generated. In this situation, we focus on the two-sided situation, and assume that we have at hand $n$ i.i.d. copies $(X_{i,j})_{i\in\Nm}\in \R^m$, $j\in \Nn$, where each $(X_{i,j})_{i\in\Nm}$ follows the location model \eqref{model-highdim}.
The $p$-values are assumed to be given by $p_i(X)=\ol{G}\left( |T(X_{i,j}, 1\leq j\leq n)|\right)$, where $T(X_{i,j}, 1\leq j\leq n)$ is some statistic, and   the (known) function $\ol{G}$ is given by  $\ol{G}(x)=\P(|T(\varepsilon_{j}, 1\leq j\leq n)| \geq x)$, $x\geq 0$, for $n$ i.i.d. copies $\varepsilon_{j}, 1\leq j\leq n$ of $\varepsilon_1$. Then, by a standard argument (see, e.g., \citealp{ABR2010b}), the joint distribution of $
(p_i(X))_{i\in\cH_0(P)}
$
can be approximated by random sign-flipping: let $\mathcal{G}=\{-1,1\}^n$ denote the group  of signs $s\in\{-1,1\}^n$ that acts on the observed $X$ in the following way:
$$
(s . X)_{i,j} = s_j X_{i,j},  \:i\in\Nm, \:j\in\Nn.
$$
Then, if $i\in\cH_0$, by symmetry, the distribution of $p_i(X)$ is equal to the one of $p_i(s.X)$, for some random sign $s$ uniformly generated in $\mathcal{G}$. As a consequence, the distribution of $
(p_i(s.X))_{i\in\cH_0(P)}$  conditionally on $X$ can act as
proxy for the distribution of $
(p_i(X))_{i\in\cH_0(P)}
$. This ``randomization property" will be formalized in detail in the next section.

Both known and unknown situations can be met in the simple Gaussian location model for which $\varepsilon\sim\mtc{N}(0,\Sigma)$ with some covariance matrix $\Sigma$ (assuming $\Sigma_{i,i}=1$ for $i\in\Nm$ for simplicity).
On the one hand, the known dependence case corresponds to the case where $\Sigma$ is known (with $\nu_m=\mathcal{N}(0,\Sigma)$). It can be met in practice in a standard Gaussian linear model or in marginal regression, see \cite{FHG2012}.
On the other hand, the unknown dependence case corresponds to the general situation where we have no information on $\Sigma$. A suitable statistics is then
$T(X_{i,j}, 1\leq j\leq n) = n^{-1/2}\sum_{j=1}^n X_{i,j}$, for which $\ol{G}(x)=2 \:\P(Z\geq x)$, $x\geq 0$, $Z\sim \mathcal{N}(0,1)$.

Also, mainly for illustrative purposes, we will use  throughout the paper the $\rho$-equi-correlated covariance matrix for which $\Sigma_{i,j}=\rho$ for $1\leq i\neq j\leq m$, for some $\rho\in[0,1]$ (either known or not).

\subsection{General framework and assumptions}\label{sec:generalassump}

Now that we have a concrete example in mind, we go 
beyond the location model by
presenting general assumptions on the $p$-value family $(p_i(X),i\in\cH_0)$. 

\paragraph{Known dependence} 
We assume that there exists a family of ``least favorable" variables $(q_i(X))_{1\leq i \leq m}$ such that for all $P \in \mathcal{P}$, 
\begin{equation}\tag{LeastFavor}\label{assumknowndep}
 \left\{\begin{array}{l}\forall i\in\cH_0(P),\:\: p_i(X) \geq q_i(X) \:\:\mbox{ $P$-a.s. }\\
\nu_m = \mathcal{D}((q_i(X))_{1\leq i \leq m}) \mbox{ does not depend on $P$.}
\end{array}\right.
\end{equation}
While \eqref{assumknowndep} is satisfied in particular in the location model (with known dependence), it encompasses some other models (e.g., scaling model).

\paragraph{Unknown dependence}
We assume that there is a finite transformation group $\mathcal{G}$ acting onto the observation set $\mathcal{X}$. 
Next, by denoting $p_{ \cH_0}(x)$ the  null $p$-value vector $(p_i(x))_{i\in \cH_0(P)}$ for  $x\in\mathcal{X}$, we assume that the joint distribution of the transformed null $p$-values is invariant under the action of any $g\in \mathcal{G}$, that is,
\begin{equation}
\forall P \in \mathcal{P},\:\: \forall g \in\mathcal{G},\:\: (p_{\cH_0}(g'.X))_{g'\in \mathcal{G}} \sim (p_{\cH_0}(g'.g.X))_{g'\in \mathcal{G}}  \label{randhypo} \tag{Rand},
\end{equation}
where $g.X$ denotes $X$ that has been transformed by $g$. This assumption has been introduced in \cite{HG2017} and is slightly weaker than the so-called randomization hypothesis of \cite{RW2005}.
It is easy to check that  \eqref{randhypo} is satisfied in the location model (with unknown dependence) for the above-mentioned sign-flipping group $\mathcal{G}=\{-1,1\}^n$, by using the symmetry of the noise. Assumption \eqref{randhypo} is also met in permutation-based two-sample multiple testing problems, as described in Section~\ref{sec:twosample}.

\section{JER control based on classical inequalities}\label{sec:basicineq}

In this section, we present an elementary approach where JER control \eqref{eq:jfwer} is derived from probabilistic inequalities that are well-known in multiple testing literature.

\subsection{Simes reference family}

\begin{theorem}[Simes and Hommel inequalities]\label{prop:simes-hommel}
  Let $(p_i(X))_{i \in \Nm}$ be a $p$-value family for the null hypotheses $(H_{0,i})_{i \in \Nm}$,
satisfying the characteristic property
  \begin{equation}\label{equ:pvalueproperty}
  \forall P \in \mtc{P}, \forall i \in \cH_0(P), \:\:\: \forall t \in [0,1], \:\P_{X\sim P}(p_i(X)\leq t)\leq t.
  \end{equation}
Then it holds that $  \forall P \in \mtc{P},$
\begin{equation}\label{equ-SimesHommel}
  \P_{X\sim P}\left( \exists k \in\{ 1,\dots,m_0\}\,:\;
  p_{(k:\cH_0)}
    \leq \frac{\alpha k}{m_0 c_m}\right) \leq \alpha\,,
  \end{equation}
where: \begin{itemize}
  \item[(i)] $c_m = C_m := \sum_{i=1}^m 1/i$ under arbitrary dependency of the $p$-value family;
  \item[(ii)] $c_m = 1$ if for all  $P\in \mtc{P}$, the $p$-value family is
  positive regression dependent on each element of the subset $\cH_0(P)$ (in short, PRDS on $\cH_0(P)$).
  \end{itemize}
  Moreover, \eqref{equ-SimesHommel} is an equality (with $c_m=1$) when the $p_i$, $i\in\cH_0(P)$, are i.i.d. $U(0,1)$.
\end{theorem}

The inequalities corresponding to items (i) and  (ii) are often referred to as the Hommel inequality \citep{Hom1983} and the Simes inequality \citep{Sim1986}, respectively.
We refer to \citet{BY2001} for a formal definition of the PRDS property. We recall that in the Gaussian model defined in Section~\ref{sec:leadingmodel} (one-sided), the PRDS assumption is valid if $\Sigma_{i,j}\geq 0$ for all $i,j\in\Nm$.

In view of \eqref{equ-JR4pvalue}, inequality \eqref{equ-SimesHommel} implies that the JER control \eqref{eq:jfwer} is satisfied for $K=m$ (under the appropriate conditions) by the reference family $\Rfam^0=(R_1^0(X),\ldots,R_m^0(X))$ given by
\begin{equation}\label{eq:ref-fam}
R^0_k(X) = \left\{i\in\Nm\::\: p_i < \frac{\alpha k}{m c_m} \right\}\,, 1\leq k \leq m.
\end{equation}
Above, we have upper-bounded $m_0$ by $m$ because $m_0$ is generally unknown.
The Hommel inequality is known to be exaggeratedly conservative, because the correction term $C_m$ is of the order of $\log(m)$. Therefore, we will only use in the sequel the reference family $\Rfam^0$ when $c_m=1$ and refer to it as the {\it Simes reference family}.
The corresponding bound is given by
\begin{align}\label{eq:VR0}
\ol{V}_{\Rfam^0}(R) &= \min_{k \in \{1, \dots, m\}}\left \{\sum_{ i \in R} \ind{p_i(X)\geq  \alpha k/m} + k-1\right\}
\:,\:\: R\subset \Nm\,.
\end{align}
This bound is considered as a baseline for our work.
As shown in Section~\ref{sec:compclosed}, this bound is in fact equivalent to
the one proposed in \cite{GS2011} for Simes local tests.

\subsection{Sharpness and conservativeness}
\label{sec:sharpness-conservativeness}

An important limitation of the reference family $\Rfam^0$ is its conservativeness and lack of adaptiveness, that is, even if $\max_{P\in \cP}\JER(\Rfam^0,P)$ is close to $\alpha$, $\JER(\Rfam^0,P)$ can be far from $\alpha$ for the $P$ that truly generated the data.
Indeed, both inequalities stated in Theorem~\ref{prop:simes-hommel} are adjusted to a \emph{worst case dependency},
thus do not adapt or take into account the dependence between the tested hypotheses.
For example, when the test statistics are strongly positively dependent, the Simes inequality may be too conservative, and the associated \iop bounds will inherit this conservativeness.

To illustrate this point, we carried out a simulation study in the Gaussian equi-correlated model where the one-sided test statistics follow the distribution $\cN(0, \Sigma)$ with $\Sigma_{ii}=1$ and $\Sigma_{ij}=\rho$ for $i \neq j$, for some $\rho\geq 0$. As noted above, this $p$-value family is PRDS. 
We consider a white setting (that is, all null hypotheses are true, $m_0=m=1,000$). In Table~\ref{tab:simes-equi}, we quantify the conservativeness of JER control in this model as the ratio of the JER actually achieved (estimated from $1,000$ simulations) to the target JER level $\alpha$ (for $\alpha=0.2$). For example, we observe that for $\rho=0.2$, the JER actually achieved by the canonical reference family $\Rfam^0$ is only 73\% of the target JER level.

\begin{table}[!htp]
\begin{tabular}{|c|c|c|c||c||c|}
      \hline
      Equi-correlation level: $\rho$ & 0 & 0.1 & 0.2 & 0.4 & 0.8\\
      \hline
      Achieved JER $\times \alpha^{-1}$ & 1.00 & 0.89 & 0.73 & 0.46 & 0.39\\
      \hline
\end{tabular}
\caption{Conservativeness of JER control based on Simes inequality in the Gaussian equi-correlated model. Here, $m_0=m=1,000$ and $\alpha=0.2$. The standard error estimate is below $0.001$ in all cases.}
\label{tab:simes-equi}
\end{table}

\subsection{Unbalancedness}\label{sec:unbal}

Let us consider a ``favorable'' case $P$ for the Simes procedure, for which the $p$-values are 
  i.i.d. uniform on $(0,1)$.
In this case, the Simes inequality is an equality
\begin{equation}\label{equ-Simesequal}
 \P_{X\sim P}\left( \exists k \in\{ 1,\dots,m\}\,:\;
  p_{(k:m)}
   < \frac{\alpha k}{m}\right) = \alpha\,,
\end{equation} 
where $p_{(k:m)}$ is the $k$-th smallest $p$-value.
In particular, the conservativeness described in Section~\ref{sec:sharpness-conservativeness} is not true here, and we might conclude that the family reference  $\Rfam^0$ given by \eqref{eq:ref-fam} can be suitably used for our aim.
However, we argue that the errors in the event described in \eqref{equ-Simesequal} are {\it not balanced} w.r.t. the parameter $k$.
As an illustration,
$
   \P(p_{(1:m)} < \alpha/m) = 1- \paren{1-\frac{\alpha}{m}}^{m} = \alpha + o(\alpha),
 $
   hence the probability of the event in \eqref{equ-Simesequal} is already almost exhausted for $k=1$. More generally, some values of the function $k\mapsto \P(p_{(k:m)} < \alpha k/m)$ are given in Table~\ref{tab:simes-unbal}  for $m=1,000$, where  $p_{(k:m)}\sim \mathrm{Beta}(k,m+1-k)$.   As a consequence, the Simes family seems to favor some of the $k$'s when controlling the JER.  In addition, the structure of this unbalancedness is somewhat arbitrary, and imposed to the user of the procedure, which may be undesirable. This phenomenon is quantified more formally in Section~\ref{append:baltemplate}, see \eqref{equ:notbal}.

  \begin{table}[!h]
    \scriptsize
    \centering
    \begin{tabular}{|c||c|c|c|c|c||c||c||c|}
      \hline
 $k$ & 1 & 2   &5  & 10 & 100\\
      \hline
 $\P(p_{(k:m)}\leq \alpha k/m)$     & $4.9 \times 10^{-2}$ & $4.7 \times 10^{-3}$   & $6.6 \times 10^{-6}$ & $1.6 \times 10^{-10}$ & $5.8 \times 10^{-93}$ \\
      \hline
    \end{tabular}
    \vspace{5mm}
    \caption{Values of  $\P(p_{(k:m)} < \alpha k/m)$ for several $k$ 
    when $p_{(k:m)}\sim \mathrm{Beta}(k,m+1-k)$,  $m=1,000$ and $\alpha=0.05$.
}
    \label{tab:simes-unbal}
  \end{table}

\section{Methodology for adaptive JER control}
\label{sec:improving-jr-control}

In this section, we aim at 
building a thresholding-based reference family $\Rfam$  for which the quantity
$
\JER(\Rfam,P)
$
is as close as possible to $\alpha$, for ``many interesting $P$s''.  To this end, we combine two approaches:

\begin{itemize}
\item incorporating the dependence structure of the noise (either known or unknown);
\item using a step-down algorithm to adapt to the unknown set $\cH_0$.
\end{itemize}

\subsection{Threshold templates}\label{sec:pivot}

We start with considering a reference family $\Rfam_\lambda$ of the form \eqref{equRkthresholdpvalues}, parametrized by $\lambda \in [0,1]$ and itself
based on a parametrized family of thresholds
$t_k(\lambda)$ which we call {\em template}.
The second step will be to to choose $\lambda=\lambda(\alpha)$ so that the JER control \eqref{eq:jfwer} is satisfied, which we call {\em $\lambda$-calibration.}

\begin{definition}\label{def:template}
A {\it one-parameter threshold template} (simply referred to as {\it template} in the sequel for short)
is a family of functions $t_k(\lambda)$, $\lambda\in[0,1]$, $1\leq k \leq K$, such that $K\in\{1,\dots,m\}$ and  for all $k\in\{1,\dots,K\}$, $t_k(0)=0$ and $t_k(\cdot)$ is non-decreasing and left-continuous on $[0,1]$. The parameter $K$ is called the {\it size} of the template.
\end{definition}

In general, a template is allowed to depend on the observation $X$.
For a given template and fixed $\lambda$, we refer to $t_k(\lambda)$, $1\leq k \leq K$, as thresholds and denote by $\Rfam_\lambda$ the associated reference family given by \eqref{equRkthresholdpvalues}. Several choices of template are possible as we will see in Section~\ref{sec:two-templates}. Here, we work with a generic, fixed template $t_k(\lambda)$, $\lambda\in[0,1]$, $1\leq k \leq K$. We denote the generalized inverse of $t_k(\cdot)$ by $t_k^{-1}(y)=\max\{x\in[ 0,1] \::\: t_k(x)\leq y\}$, for any $y\in \R\cup\{-\infty,+\infty\}$.

Since $t_k(\cdot)$ is monotonic, for any $p$-value family $\{p_i, i\in \Nm\}$, we have $t_k(\lambda)> p_{(k:\cH_0)}$ if and only if $\lambda> t_k^{-1}(p_{(k:\cH_0)})$. Hence, in view of \eqref{equ-JR4pvalue}, we obtain
\begin{align}
\JER(\Rfam_\lambda,P)&=
\P_{X\sim P}\bigg(\exists k \in \{ 1,\dots,K\wedge m_0\} \::\: p_{(k:\cH_0)} < t_k(\lambda) \bigg)
\nonumber\\
&=\P_{X\sim P}\bigg(\exists k \in \{ 1,\dots,K\wedge m_0\} \::\: t_k^{-1}\left(p_{(k:\cH_0)}\right) < \lambda \bigg).
\nonumber
\end{align}
This proves the following result.

\begin{lemma}\label{lem:pivotal}
Consider a general $p$-value model and any (possibly data-dependent) template $t_k(\lambda)$, $\lambda\in[0,1]$, $1\leq k \leq K$.
Then, for any $\lambda\in[0,1]$, the error rate \eqref{eq:jr} of the reference family $\Rfam_\lambda$
given by \eqref{equRkthresholdpvalues}
 can be written as follows: for any $P\in\cP$, \begin{align}
\JER(\Rfam_\lambda,P)&=\P_{X\sim P}\bigg(\min_{1\leq k \leq K\wedge m_0}\left\{ t_k^{-1}\left(p_{(k:\cH_0)}(X)\right)\right\} < \lambda \bigg) \,.\label{eq:pivotal}
\end{align}
\end{lemma}

\subsection{Single-step and step-down procedures by $\lambda$-calibration}\label{sec:lambdaadj}

The  JER control~\eqref{eq:jfwer} can now be achieved by choosing $\lambda$  in an appropriate way.
\begin{definition}\label{def:lambdaadjust}
Given a threshold template $t_k(\lambda)$, $\lambda\in[0,1]$, $1\leq k \leq K$,
a (possibly data-dependent) functional $\lambda(\alpha,A)$, $\alpha\in(0,1), A\subset \Nm$, is called a $\lambda$-calibration if it is non-increasing in $A$, that is,
\begin{equation}\label{NI}
\forall \alpha\in(0,1), \forall A,A'\subset \{1,\dots,m\}, \mbox{ with } A\subset A', \:\:\:  \lambda(\alpha,A')\leq \lambda(\alpha,A),
\end{equation}
and  satisfies $\forall \alpha\in(0,1)$, $\forall P \in \cP$,
\begin{equation}\label{lambdacalibrcontrol}
 \P_{X\sim P}\bigg(\min_{1\leq k \leq K\wedge m_0}\left\{ t_k^{-1}\left(p_{(k:\cH_0(P))}(X)\right)\right\} < \lambda(\alpha,\cH_0(P)) \bigg) \leq \alpha.
 \end{equation}
\end{definition}

Two examples of possible $\lambda$-calibrations will be provided in Sections~\ref{sec:lambdaadjknown} and~\ref{sec:lambdaadjunknown}. In  the remaining of this section, we consider that some $\lambda$-calibration is given.

The dependence of the calibration on the set $A$ adds extra flexibility which will allow us to
apply a step-down principle and get a more accurate procedure. 
A consequence of Lemma~\ref{lem:pivotal} is that
the procedure $\Rfam_{\lambda(\alpha,\cH_0)}$ has a controlled JER (given a template and a calibration),
in other words taking $A=\cH_0$ provides an ``oracle'' calibration, but
since $\cH_0$ is unknown, $\lambda(\alpha,\cH_0)$ cannot be used.
However, a consequence of \eqref{NI} is that $\lambda(\alpha,\Nm)\leq \lambda(\alpha,\cH_0)$, so that $\lambda(\alpha,\Nm)$ can be used as a (single-step) conservative substitute for $\lambda(\alpha,\cH_0)$. This provides the following result.

\begin{proposition}\label{prop:JRadjust}
In the framework of Lemma~\ref{lem:pivotal}, consider $\lambda(\alpha)=\lambda(\alpha,\Nm)$ for some $\lambda$-calibration as in Definition~\ref{def:lambdaadjust}. Then the procedure  $\Rfam_{\lambda(\alpha)}$ controls the JER criterion at level $\alpha$ in the sense of \eqref{eq:jfwer}.
\end{proposition}

Above, the fact that $\lambda (\alpha,\Nm)$ is smaller than $\lambda(\alpha,\cH_0)$ induces a loss in the JER control. This loss can sometimes by substantial, as illustrated with numerical experiments in Section~\ref{sec:num}; this effect is further studied theoretically in Section~\ref{sec:effects-step-down}.
This loss can be reduced by using $\lambda(\alpha,\wh{A})$, where $\wh{A}$ is the output of the the following step-down algorithm.

\IncMargin{1em}
\begin{algorithm}[H]
  \label{algoSD}
  $j \leftarrow 0$ \;
  $A^{(0)} \leftarrow \Nm$\;
  \Repeat{$A^{(j)} =  A^{(j-1)}$}{
    $j \leftarrow j+1$ \;
    $\lambda_j  \leftarrow \lambda(\alpha,A^{(j-1)})$ \;
    $A^{(j)}  \leftarrow  \left\{i\in\Nm \::\: p_i(X)\geq  t_1(\lambda_j)\right\}$ \;
  }
  \Return{$A^{(j)}$}\;
  \caption{General step-down algorithm}
\end{algorithm}

While the update of $A^{(j)}$ only depends on $t_1(\cdot)$ in Algorithm \ref{algoSD}, $\wh{A}$ may depend on all the $t_k$'s through the functional $\lambda(\alpha,\cdot)$. The following result is proved in Section~\ref{proof:improving-jr-control}.

\begin{proposition}\label{prop:JRadjustSD}
In the framework of Lemma~\ref{lem:pivotal}, consider any $\lambda$-calibration as in Definition~\ref{def:lambdaadjust} and compute $\wh{A}$ by Algorithm~\ref{algoSD}. Then the procedure  $\Rfam_{\lambda(\alpha,\wh{A})}$ controls the JER at level $\alpha$ in the sense of \eqref{eq:jfwer}.
\end{proposition}

\begin{remark}\label{rem:RW}
 When we choose $K=1$, Algorithm~\ref{algoSD} reduces to the usual FWER controlling step-down algorithm (see, e.g., \citealp{RW2005}).
\end{remark}

\subsection{Valid $\lambda$-calibration for known dependence}\label{sec:lambdaadjknown}

Let us focus on the situation where the dependence is known, see  Section~\ref{sec:generalassump}. The template is assumed to be deterministic in this section. Assumption \eqref{assumknowndep} and 
 Lemma~\ref{lem:pivotal} thus give
\begin{align} 
\label{eq:JER-upper-bound-known-dep}
\JER(\Rfam_\lambda,P)&\leq
\P_{q \sim \nu_m}\bigg(
\min_{1\leq k \leq K\wedge m_0}\left\{ t_k^{-1}\left(q_{(k:\cH_0)}\right)\right\} < \lambda \bigg) \,,
\end{align}
which provides the following valid $\lambda$-calibration: for all $A\subset\{1,\dots,m\}$,
\begin{align}
\lambda(\alpha,A) &=\max\bigg\{ \lambda \geq 0\telque \P_{q\sim \nu_m}\bigg(\min_{1\leq k \leq K\wedge |A|}\left\{ t_k^{-1}\left(q_{(k:A)}\right)\right\} < \lambda \bigg) \leq \alpha \bigg\}.\label{lambdaadj}
\end{align}
Property~\eqref{NI} can be easily checked.
Note that $\lambda(\alpha,\cdot)$ depends on $\nu_m$ and on the template, although it is not explicit from the notation for short.
We have proved the following result.
\begin{theorem}[$\lambda$-calibration for known dependence]\label{prop:knowndep}
Consider any $p$-value family  satisfying \eqref{assumknowndep}, a deterministic template and the associated reference family $\Rfam_\lambda$.
Then the  (deterministic) functional $\lambda(\cdot,\cdot)$ defined by \eqref{lambdaadj} is a $\lambda$-calibration in the sense of Definition~\ref{def:lambdaadjust}
and thus $\Rfam_{\lambda(\alpha,\Nm)}$ and $\Rfam_{\lambda(\alpha,\wh{A})}$ both control the JER at level $\alpha$. 
\end{theorem}

\subsection{Valid $\lambda$-calibration for unknown dependence}\label{sec:lambdaadjunknown}

Let us consider now the case where the dependence is unknown, see  Section~\ref{sec:generalassump}.
The template is still assumed to be deterministic in this section.
We use the notation defined therein and in particular assumption  \eqref{randhypo}.
 Let us consider a (random) $B-$tuple $(g_1,g_2,\dots,g_B)$ of $\mathcal{G}$ (for some $B\geq 2$), where $g_1$ is the identity element of $\mathcal{G}$ and $g_2,\dots,g_B$ have been drawn  (independently of the other variables) as i.i.d. variables, each being uniformly distributed on $\mathcal{G}$.

Let us consider some template $t_k(\cdot)$, $1\leq k \leq K$, and, for short,   denote for all $A\subset\Nm$,
$$
\Psi(X,A)= \min_{1\leq k \leq K\wedge |A|}\left\{ t_k^{-1}\left(p_{(k:A)}(X)\right)\right\}.
$$
Now introduce the (data-dependent) $\lambda$-calibration
\begin{align}
\lambda(\alpha, A) &=\max\bigg\{ \lambda \geq 0\telque B^{-1}\sum_{j=1}^B \ind{ \Psi(g_j.X,A) < \lambda} \leq \alpha \bigg\}.\label{lambdaadjrand}
\end{align}
In practice, we can compute this functional easily as $\lambda(\alpha, A)=\Psi_{(\lfloor \alpha B\rfloor+1)}$
where $\Psi_{(1)}\leq \Psi_{(2)}\leq \cdots \leq \Psi_{(B)}$ denote the ordered sample $(\Psi(g_j.X,A), 1\leq j \leq B)$.
Then the following result holds and is proved in Section~\ref{sec:proof:prop:unknowndep}.

\begin{theorem}[$\lambda$-calibration for unknown dependence]\label{prop:unknowndep}
  \label{PROP:UNKNOWNDEP}
Consider any $p$-value family  satisfying \eqref{randhypo}, a deterministic template and the associated reference family $\Rfam_\lambda$. 
Then the (data-dependent) functional $\lambda(\cdot,\cdot)$ defined by \eqref{lambdaadjrand} is a $\lambda$-calibration in the sense of Definition~\ref{def:lambdaadjust}
and $\Rfam_{\lambda(\alpha,\Nm)}$ and $\Rfam_{\lambda(\alpha,\wh{A})}$ both control the JER at level $\alpha$
\end{theorem}

A related idea has been proposed independently by \cite{HGS2017} to build confidence envelopes for the False Discovery Proportion.

\section{Application : two examples of template-based reference families}
\label{sec:two-templates}

In this section, we apply the methodology presented in the previous section for two particular instances of templates. Throughout this section, the $\lambda$-calibration functional $\lambda(\alpha,A)$ is either given by \eqref{lambdaadj} (known dependence) or by \eqref{lambdaadjrand} (unknown dependence).

\subsection{Linear template}\label{sec:linreffam}

We define the {\it linear template} (of size $K$) by
\begin{equation}\label{equ:simesfam}
t^{L}_k(\lambda) = \lambda k/m, \:\:\:\lambda\in[0,1], \:\:\:1\leq k \leq K.
\end{equation}
Hence we have $(t^{L}_k)^{-1}(u)=1\wedge (\frac{m}{k}u)$ which corresponds to a specific $\lambda$-calibration denoted by $\lambda^{L}(\alpha,A)$.
For each $K$, this gives rise to two new reference families:
\begin{itemize}
\item[\textbullet] The {\it single-step linear reference family} (of size $K$), denoted $\Rfam^{L}$, is given by
 $\Rfam^L=(R_1^L(X),$ $\ldots,R_K^L(X))$, where
 \begin{equation}\label{eq:1SL-ref-fam}
R^L_k(X) = \left\{i\in\Nm\::\: p_i < \lambda^{L}(\alpha,\Nm) \frac{ k}{m} \right\}\,, 1\leq k \leq K.
\end{equation}
\item[\textbullet] The {\it step-down linear reference family} (of size $K$), denoted $\Rfam^{L, sd}$, is given by
 $\Rfam^{L,sd}=$ $(R_1^{L,sd}(X),\ldots,$ $R_K^{L,sd}(X))$, where
\begin{equation}\label{eq:1SDL-ref-fam}
R^{L,sd}_k(X) = \left\{i\in\Nm\::\: p_i < \lambda^{L}(\alpha,\hat{A}) \frac{ k}{m} \right\}\,, 1\leq k \leq K,
\end{equation}
where $\hat{A}$ is derived from Algorithm~\ref{algoSD}, used with $\lambda(\cdot)=\lambda^L(\cdot)$ and $t_1(\cdot) = t^{L}_1(\cdot)$.
\end{itemize}

Theorems~\ref{prop:knowndep}~and~\ref{prop:unknowndep} ensure that the reference families $\Rfam^{L}$ and $\Rfam^{L, sd}$ both control the JER at level $\alpha$ both in the known and unknown dependent case. 
\magen{The magnitude of $\lambda^{L}(\alpha,\Nm)$ is studied in Section~\ref{sec:magnitude} in a simple case. It shows that our $\lambda$-calibration adapts to the dependence structure and addresses the conservativeness issue raised in Section~\ref{sec:sharpness-conservativeness}.}

\subsection{Balanced template}\label{sec:balreffam}

Considering a linear template is not always appropriate: as mentioned above, under independence and $K=m$,  $\Rfam^{L}$ corresponds to the Simes reference family $\Rfam^{0}$ \eqref{eq:ref-fam}, and thus suffers from a kind of unbalancedness, as underlined in Section~\ref{sec:unbal}.
\magen{Ideally, a {\em balanced} reference family $R_k$ would have the property that
  $\P(|R_k| \geq k)$ is a constant not depending on $k=1,\ldots,K$. While strict balancedness
  seems out of reach, since these probabilities depend on $\cH_0$, we can ensure balancedness under the full null configuration ($\Nm=\cH_0$) by calibrating the template
  as a quantile at a common level for all $k$, as follows.}
For each $k\in\Nm$, let us define
$$
\left\{\begin{array}{ll} F_k(x)=\P_{q\sim \nu_m}(q_{(k:m)}\leq  x)& \mbox{ (known dep.)}\\
F_k(x)=B^{-1}\sum_{j=1}^B \ind{ p_{(k:m)}(g_j.X)\leq x}& \mbox{ (unknown dep.)}
\end{array}\right., \:\:x\in [0,1].
$$
The {\it balanced template} (of size $K$) is then given by
\begin{equation}\label{equ:skbal}
t^B_k(\lambda) = F_k^{-1}(\lambda) = \min\{ x \in [0, 1] \::\: F_k(x)\geq \lambda\}, \:\: \mbox{ with }\:\: k\in\{1,\dots,K\}.
\end{equation}
From an intuitive point of view, for each $k$, the threshold $t^B_k(\lambda)$ corresponds to a procedure controlling the  $\kFWER$ at level $\lambda$. It is straightforward to check that $t^B_k(\cdot)$ fulfills the requirements of Definition~\ref{def:template} while $(t^B_k)^{-1}(x)=F_k(x)$ for all $x\in [0,1]$. This corresponds to a specific $\lambda$-calibration denoted by $\lambda^{B}(\alpha,A)$.
For each $K$, this gives rise to two new reference families:
\begin{itemize}
\item[\textbullet] The {\it single-step balanced reference family} (of size $K$), denoted $\Rfam^{B}$, is given by
 $\Rfam^B=(R_1^B(X),$ $\ldots,R_K^B(X))$, where
 \begin{equation}\label{eq:1SB-ref-fam}
R^B_k(X) = \left\{i\in\Nm\::\: p_i < t_k^B(\lambda^{B}(\alpha,\Nm))  \right\}\,, 1\leq k \leq K.
\end{equation}
\item[\textbullet] The {\it step-down balanced reference family} (of size $K$), denoted $\Rfam^{B, sd}$, is given by
 $\Rfam^{B,sd}=$ $(R_1^{B,sd}(X),$ $\ldots,R_K^{B,sd}(X))$, where
\begin{equation}\label{eq:1SDB-ref-fam}
R^{B,sd}_k(X) = \left\{i\in\Nm\::\: p_i < t_k^B(\lambda^{B}(\alpha,\hat{A})) \right\}\,, 1\leq k \leq K,
\end{equation}
where $\hat{A}$ is derived from Algorithm~\ref{algoSD}, used with $\lambda(\cdot)=\lambda^B(\cdot)$ and $t_1(\cdot) = t^{B}_1(\cdot)$.
\end{itemize}
We give in section Section~\ref{MCapprox} a detailed construction of the reference families $\Rfam^{B}$ and $\Rfam^{B,sd}$. Theorem~\ref{prop:knowndep} ensures that both of these reference families control the JER at level $\alpha$ in the case of a known dependence. 

However, for unknown dependence, Theorem~\ref{prop:unknowndep} cannot be directly applied to the balanced template. Indeed, although this is not acknowledged by the notation for simplicity, $F_k$ and thus $t^B_k(\lambda)$ depend on the observation $X$. Our proof does not generalize easily to such a data-dependent rejection template, although the numerical experiments of Section~\ref{sec:num} suggest that the JER control is also valid in that situation. 

\begin{remark} 
The step-down refinement can be substantial for a balanced template, as further discussed in Section~\ref{sec:effects-step-down}.
\end{remark}

\begin{remark} 
By considering the two-sample setting with unknown dependency structure (see Section~\ref{sec:twosample}) our balanced procedure  is related to the work of \cite{Mein2006}, where permutations are used to build FDP confidence envelopes. However, there appears to be a gap in the theoretical analysis justifying the validity of such an approach
(Theorem~1 of \citealp{Mein2006}, more specifically Equation~(12) there), which seems to have been overlooked so far. 
The reason is similar to the one making our proof not cover the case of a data-dependent template  
  $t_k(X,\lambda)$: the fact that
for all $\lambda$ and ${g}\in\mtc{G}\,$, $(t_k(g.X,\lambda))_{1\leq k \leq K}=(t_k(X,\lambda))_{1\leq k \leq K}$ and $(p_{i}(g.X))_{i\in\cH_0} \sim (p_{i}(X))_{i\in\cH_0} $,
does not imply (in general) equality of the joint distributions $((t_k(X,\lambda))_{1\leq k \leq K},(p_{i}(X))_{i\in\cH_0})$ and $((t_k(g.X,\lambda))_{1\leq k \leq K},(p_{i}(g.X))_{i\in\cH_0})$. 
\end{remark}

\section{Numerical experiments}\label{sec:num}

We report numerical experiments performed in the two-sided location model~\eqref{model-highdim} described in Section~\ref{sec:leadingmodel} in the case of an \emph{unknown dependence}. The observations $(X_{i,j})_{i\in\Nm}\in \R^m$, $j\in \Nn$ are distributed as $\rho$-equi-correlated, and the test statistics for $i \in \Nm$ is $T(X_{i,j}, 1\leq j\leq n) = n^{-1/2}\sum_{j=1}^n X_{i,j}$. We use sign-flipping (as described in that section) to approximate the joint distribution of the test statistics under the null. The location parameter is set to $\mu_i=n^{-1/2} \ol{\mu}\: \ind{i \in \cH_1}$, where $\ol{\mu} >0 $ quantifies the signal-to-noise ratio (SNR). We have also performed experiments in the same model but assuming \emph{known dependence}, in order to illustrate Theorem~\ref{prop:knowndep}. The results of these experiments are quite similar to those reported here for unknown dependence. 

\subsection{JER control}
\begin{figure}[!htp]
  \centering
\includegraphics[width=0.99\textwidth]{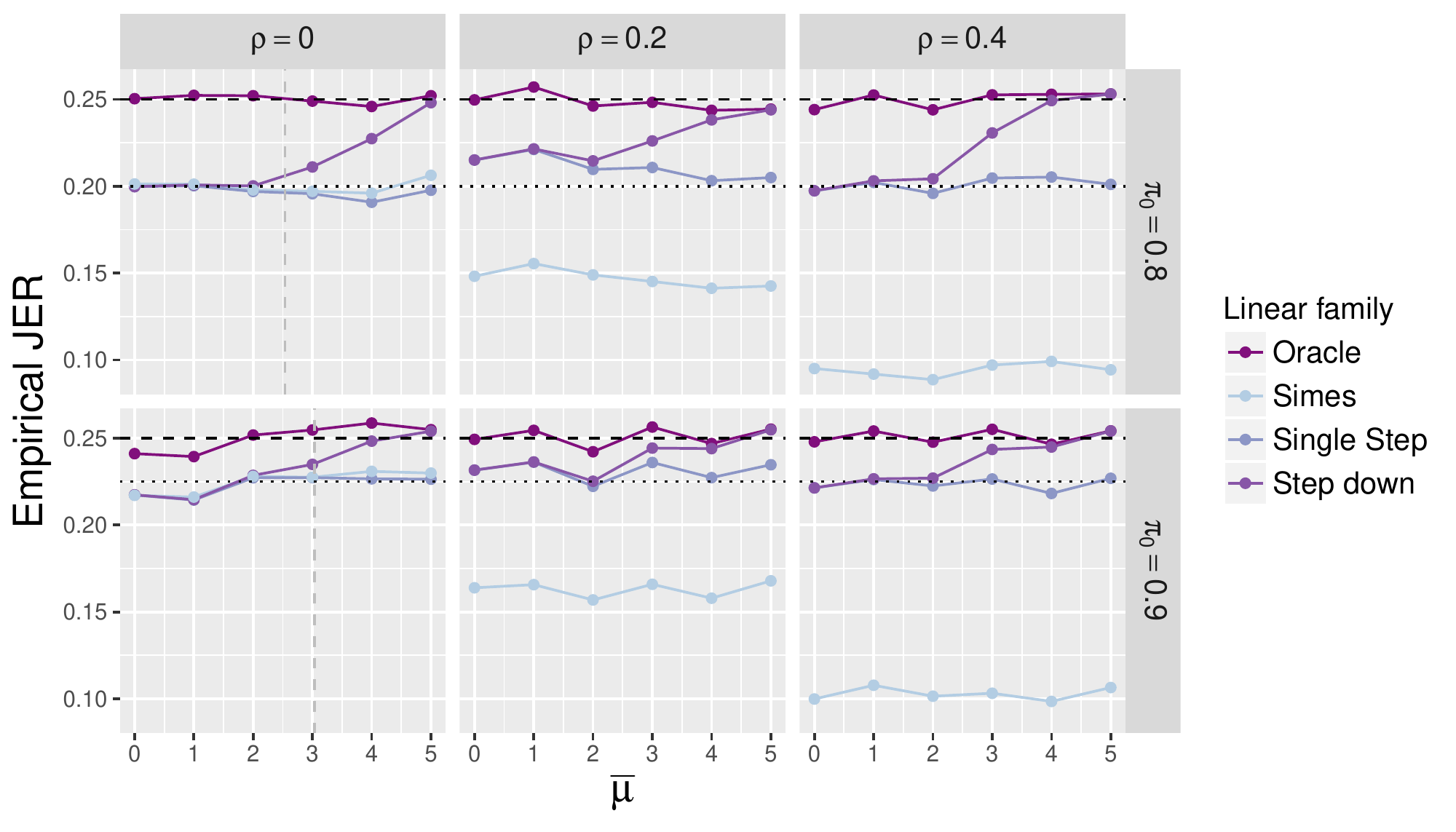}
  \caption{JER control based on the linear template for equi-correlated test statistics.}
  \label{fig:simes-equi-flip} 
\end{figure}

\begin{figure}[h!]
  \centering
\includegraphics[width=0.99\textwidth]{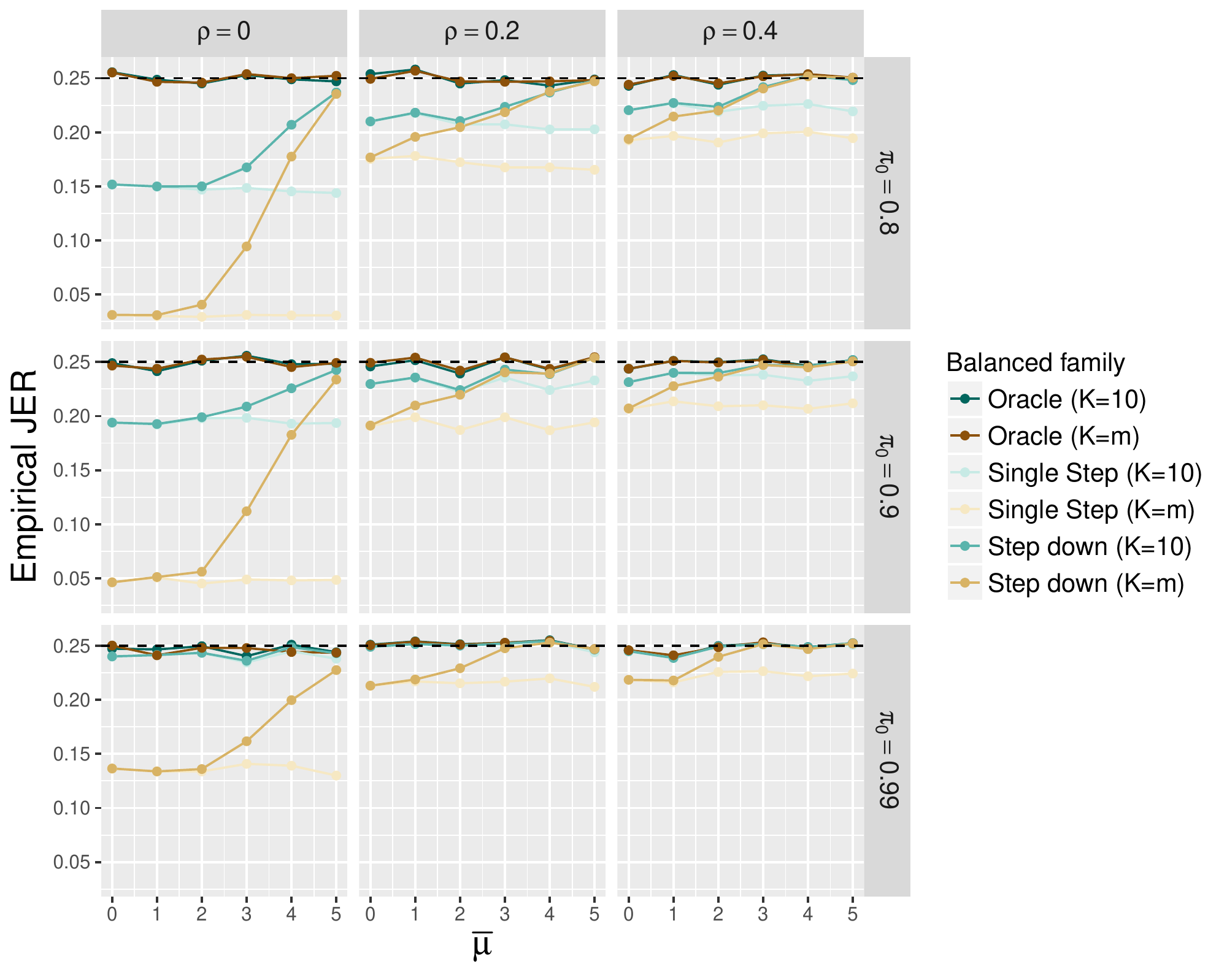}
  \caption{JER control based on the balanced template for equi-correlated test statistics, with $K=m$ and $K=10$.}
  \label{fig:balanced-equi-flip}
\end{figure}

The target JER level is set to $\alpha=0.25$, and the simulation parameters are: $m=n=1,000$, $\rho \in \{0, 0.2, 0.4\}$, $\pi_0 \in \{0.8, 0.9, 0.99\}$ (corresponding to $m_1 \in \{200, 100, 10\}$), and $\ol{\mu} \in \{0, 1, 2, 3, 4, 5\}$. For each setting, we report the empirical JER achieved, that is, the proportion of simulation runs (out of a total of $10,000$ runs) for which $|R_k(X) \cap \cH_0(P)| > k$ for at least one $k \in \{1, \dots, K\}$. 
The results are summarized by Figure~\ref{fig:simes-equi-flip} for the linear template, and by Figure~\ref{fig:balanced-equi-flip} for the balanced template. Each figure is a matrix of panels, where each row corresponds to one value of the sparsity parameter $\pi_0$, and each column corresponds to one value of the equi-correlation parameter $\rho$. In each panel, the empirical JER achieved by several procedures is displayed as a function of the signal-to-noise ratio parameter $\ol{\mu}$.
The target JER level $\alpha$ is represented by a horizontal dashed line, and for the linear template, the level $\pi_0 \alpha$ is represented by a horizontal dotted line. In both figures, each color corresponds to a different $\lambda$-calibration:
\begin{center}
  \begin{tabular}{ccc}
    single-step & Step down & Oracle\\ \hline
    $\lambda(\alpha,\Nm)$ & $\lambda(\alpha,\wh{A})$ & $\lambda(\alpha,\cH_0)$
  \end{tabular}
\end{center}
Additionally, for the linear template, ``Simes'' corresponds to $\lambda=\alpha$ (no $\lambda$-calibration).
Figure~\ref{fig:simes-equi-flip} illustrates that the JER is controlled at the target level $\alpha$ in all situations for the linear template, which is expected according to Proposition \ref{prop:unknowndep}. Oracle calibration yields exact JER control, up to sampling fluctuations.  As discussed in Section \ref{sec:sharpness-conservativeness}, the Simes reference family with parameter $\alpha$ yields JER equal to $\pi_0 \alpha$ under independence ($\rho=0$), while it is more conservative under positive dependence $\rho>0$. Single-step $\lambda$-calibration addresses this conservativeness by adapting to the (unknown) dependence: it yields JER control at $\pi_0 \alpha$ in all settings considered. Finally, as the signal-to-noise ratio $\ol{\mu}$ gets larger, the step-down $\lambda$-calibration yields a JER closer to the nominal level $\alpha$ in non-sparse situations ($\pi_0 \in \{0.8, 0.9\})$. In a sparse situation ($\pi_0 = 0.99)$, corresponding to $m_1 =10$, the single-step procedure is already quite sharp and essentially indistinguishable from its Oracle counterpart, so we decided to omit this setting from Figure~\ref{fig:simes-equi-flip}.

The results for the balanced template are summarized by Figure~\ref{fig:balanced-equi-flip}. First, the JER is empirically controlled at the target level $\alpha$ in all situations. This is worth noting because as discussed in the preceding section, our results do not cover the case of unknown dependence for the balanced template. Looking at the (brown) curves corresponding to $K=m$, single-step $\lambda$-calibration leads to a much more conservative JER control than for the linear template, especially under independence or for small values of $\rho$, even when $\pi_0$ is close to one. For example, when $\pi_0=0.99$ ($m_1=10$ out of $m=1,000$), the JER achieved by the single-step $\lambda$-calibration of the balanced family is of the order of $\alpha/2 (\ll \pi_0 \alpha)$. When the signal-to-noise ratio is large, our proposed step-down adjustment catches up with the target JER level. This effect is further discussed and formalized in Section~\ref{sec:effects-step-down}. 

Interestingly, the JER control offered by the balanced family with $K=10$ (green curves in Figure~\ref{fig:balanced-equi-flip}) is much less conservative than with $K=m$, even for the single-step $\lambda$-calibration. The magnitude of the $\lambda$-adjustment is further discussed in Section \ref{sec:magnitude}, and the question of how to choose $K$ is discussed in Section~\ref{se:literature}.

\paragraph*{Additional numerical experiments} The experiments reported here are carried out only in the equi-correlated setting and assuming that the mean signal under the alternative is constant: $\mu_i=\ol{\mu}$ for all $i \in \cH_1$. We have performed other experiments, where $\mu_i$ is uniformly distributed between $0$ and $\ol{\mu}$, and/or where the test statistics have a Toeplitz covariance, for which $\Sigma_{i,j}=|i-j|^{\theta}$, where $\theta \in \{-2, -1, -0.5, -0.2\}$ controls the range of dependency. The results obtained for both types of signals and for both types of dependency are qualitatively similar, so we have only reported the results for the parameter combination: constant signal/equi-correlated dependency.

\subsection{Power}
\label{sec:power}\label{sec:averaged-power}

In the preceding section, the quality of a JER controlling procedure is quantified by the tightness of its JER control. We now compare some JER controlling procedures in terms of power. This comparison is made under independence for simplicity. We focus on the step-down linear reference family \eqref{eq:1SDL-ref-fam} with $K=m$, and the step-down balanced  reference family \eqref{eq:1SDB-ref-fam} with $K \in \{10, 2m_1, m\}$.
We consider a notion of power, referred to as ``averaged power'', that takes into account the amplitude of the lower bound $\ol{S}_\Rfam( \cdot)$.
Let us define for some selected set $R\subset \Nm$ (possibly data dependent),
\begin{equation}\label{equ:Pow2}
\Pow(\Rfam,P,R)=\E\left(\frac{\ol{S}_\Rfam( R)}{|R\cap \cH_1(P)|}\:\bigg|\: |R\cap \cH_1(P)|>0\right)\,,
\end{equation}
where we recall that $\ol{S}_\Rfam(R)= |R| - \ol{V}_\Rfam(R)$.
The following selected sets  $R\subset \Nm$ are considered:
\begin{itemize}
\item[(a)] $R=\Nm$. In this case, the averaged power $\Pow(\Rfam,P,R)$ measures the (relative) performance of $\ol{S}_\Rfam( \Nm)$ as an estimator of $m_1(P)=|\cH_1(P)|$;
\item[(b)] $R_0=\{i \in \Nm: p_i \leq 0.05\}$, and $R$ is a random selection of half of the items of $R_0$. Each hypothesis is given a selection probability proportional to the rank of its $p$-value;
\item[(c)] Same as (b) with $R_0$ corresponding to the rejections of the BH procedure at level $0.05$.
\end{itemize}
In (b)-(c) above, the sets $R$ are thought to be typical possible choices for the user. We chose to give non-uniform selection probabilities in order to favor sets enriched in lower $p$-values.
The parameter $\pi_0$ is taken in the range $\pi_0 \in \{0.8, 0.9, 0.99\}$. We set $\ol{\mu} = \sqrt{-4 \log (1-\pi_0)}$ in order to specifically focus on situations where the signal strength lies just above the estimation boundary, which would correspond to $\ol{\mu} = \sqrt{-2 \log (1-\pi_0)}$, see \cite{DJ2004}.

\begin{figure}[h!]
  \centering
\includegraphics[width=0.99\textwidth]{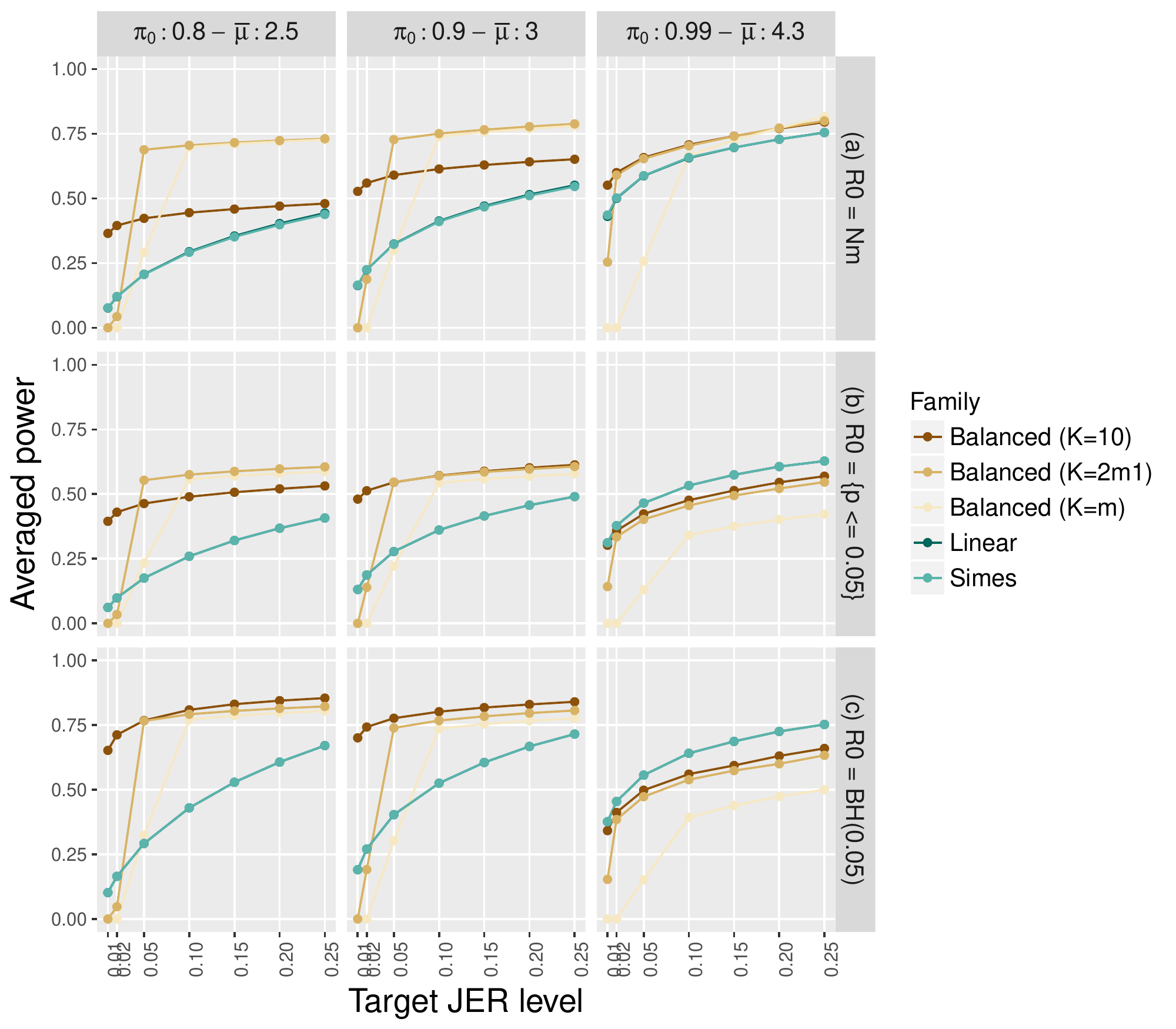}
\caption{Averaged power of JER controlling procedures for independent test statistics.}
   \label{fig:averaged-power}
\end{figure}

The results are displayed in Figure~\ref{fig:averaged-power}. The average power of the Simes family (light green) and of the reference families obtained by single-step and step-down $\lambda$-calibration of the linear template (dark green) are almost identical.  This is consistent with the results displayed in the first column of Figure~\ref{fig:simes-equi-flip}, where the three families achieve very similar JER levels for $\ol{\mu} \leq \sqrt{-4 \log (1-\pi_0)}$; this value of  $\ol{\mu}$ is shown by a dashed gray vertical line. 
Overall, the averaged power obtained from the balanced template is substantially larger than the averaged power obtained from the linear template. While neither template uniformly dominates the other one, the only situation where the linear template is more powerful is under the most sparse scenario ($\pi_0=0.99$), for the two user-defined rejection sets (b) and (c).
 In particular, the first row of panels in Figure~\ref{fig:averaged-power} indicates that, except for a very low target JER ($\alpha \leq 0.02)$, the bound $\ol{S}_\Rfam( \Nm)$ obtained from the balanced template provides a better estimator of $m_1(P)=|\cH_1(P)|$ than the linear template. These experiments also show that, as expected, the choice of $K$ can improve the performance of the balanced procedure. Some suggestions for choosing $K$ are discussed below.

\newpage
\section{Discussion}
\label{se:literature}

\subsection{Choosing the size $K$}
\label{disc:K}
While the choice $K=m$ seems a priori natural, we have shown throughout this paper that it induces some conservativeness (via the $\lambda$-calibration): choosing a smaller value for $K$ can yield a tighter post hoc bound.
This effect is particularly marked in the case of the balanced template when $p$-values are close to independent (see Figure~\ref{fig:balanced-equi-flip}). The choice of $K$ is therefore quite important in practice. We underline the following plausible scenarios:
\begin{itemize}
\item if the user has an {\it a priori maximum amount of tolerated false discoveries}, then $K$ can be set taken equal to that value. This comes from the following fact:
{\it let $K_0\in\mathbb{N}$ and assume $\Rfam=(R_i(X))_{1 \leq i \leq K}$ is a reference family  (using $\zeta_i=i-1$) satisfying JER control.
Consider any set $R\subset \Nm$ such that $\ol{V}_{\Rfam}( R) \leq K_0 < K$. Then we have $\ol{V}_\Rfam( R)=\ol{V}_{\Rfam^{(K_0)}}( R)$, where $\Rfam^{(K_0)} = (R_i(X))_{1 \leq i \leq K_0+1}$.}
In words, if the user is only interested in rejected sets $R$ where the bound on the number of false positives is less than $K_0$, then the family size $K$ can safely be taken equal to $K_0+1$.
\item if the user has some upper bound $\ol{m}_1$ on the number of false hypotheses as prior information, it seems reasonable to take $K_0 = \ol{m}_1$ above (a larger number of false discoveries would mean that more than 50\% of the hypotheses in the rejected set are false discoveries). The case $K=2m_1$ considered in our numerical experiments can be interpreted as such a scenario (assuming a known prior rough upper bound $\ol{m}_1=2m_1$).
\end{itemize}
Designing a theoretically founded data-dependent choice of $K$ is an interesting direction for future efforts.
Let us also mention that an alternative direction to the choice of $K$ is to introduce some smooth decay in the violation probability $\P(|R_k|\geq k)$ as $k$ grows.

\subsection{Step-down algorithm}

The principle of the step-down Algorithm~\ref{algoSD} is to approach the oracle value $\lambda(\alpha,\cH_0)$
by iterative approximations $\lambda(\alpha,\wh{A})$. Here the template $t_k(\cdot)$ is fixed once for all. A seemingly natural extension is to allow the template $t_k(\cdot,A)$ to also depend on subsets $A\subset \Nm$ and to apply the step-down algorithm to the template as well as $\lambda$, that is, consider  at each step
$t_k(\cdot,\wh{A})$, then apply the $\lambda$-calibration step.
For instance, for the balanced rejection template, one could define $t^B_k(\lambda,A)$
as the $\lambda$-quantile of $q_{k:A}$.
From a theoretical point of view however, it turns out that the corresponding combined threshold (depending on $\cH_0$ both through $t_k$ and $\lambda$) loses the monotonicity property with respect to $\cH_0$. Hence, our current proof does not extend to that situation and we do not know if  the corresponding JER is controlled at level $\alpha$. This is an interesting (but challenging) issue.

\subsection{Choice of the reference family}
\label{disc:choice-ref-ram}
\magen{In the general setting presented in Section~\ref{se:JRglobal},
although the aim is to obtain a uniform guarantee for any possible rejected set,
a tradeoff is implicitly present in the
choice of the reference family. The post hoc bounds \eqref{eq:optbound1}, \eqref{eq:aug}
can be understood as interpolation bounds relating an arbitrary $R$ to sets of the reference family $\Rfam$, so that generally speaking they will be more accurate for rejection sets that
are ``well approximated'' by sets of the reference family. From the definition of the JER
control \eqref{eq:jfwer}, it is clear that there is a tradeoff between the cardinality of the reference family and the conservativeness of the bound, which requires a uniform control over the family.
Depending on the specific application, reference families corresponding to
different expected tradeoffs 
can be considered. In the running example considered in this paper,
the choice of $K$ (discussed above) represents precisely such a tradeoff; so does the choice of the calibration function, as we have already argued. Adequate choice of reference families for
specific applications and goals, and an appropriate notion of which sets well approximated by the reference family, remains an important avenue to explore.
}

\subsection{Principled use of user-agnostic bounds and admissible sets}
\label{disc:admissibility}

\magen{This point stems from an
insightful remark by an anonymous reviewer.
If there are no constraints on the rejected set $R$ selected by the user,
and a post hoc bound $V(\cdot)$ is available, it seems sensible to require
that one should not be able to add hypotheses to
the rejected set without increase of the bound on false discoveries, nor exclude
hypotheses from it without decrease of the bound on true discoveries; otherwise the choice of $R$ would obviously be suboptimal given
the information given by the bound. Formally, call $R$ admissible with respect to bound $V(\cdot)$ if
\begin{itemize}
\item[(i)] $\forall R' \supsetneq R,$ $V(R')>V(R)$;
\item[(ii)] $\forall R' \subsetneq R,$ $S(R')<S(R)$.
\end{itemize}
We leave to the reader to check the following result: {\em
  the only sets admissible with respect to $\ol{V}_\Rfam$ (of \eqref{eq:aug}) belong to
  the reference family.} (In particular, for nested reference families,
only the reference sets are admissible with respect to the optimal post hoc bound
$V^*_\Rfam$). This property emphasizes the role played by the choice of reference
family --- while also putting into question to allow 
rejection sets not belonging to it in the first place. 
Concerning this last point, we argue that additional constraints (sometimes only implicitly
defined by the selection procedure used) often restrict the rejection sets under consideration of the user
(this is the case in the two exemplary applications mentioned in the introduction).
In such a  situation, the reference sets might not satisfy the constraints, which
justifies the interest of a bound for more general $R$s.
One may in this case adapt the above definition of admissible sets by restricting
comparisons to sets satisfying the constraints;  which sets are then
admissible would have to be investigated in specific situations.}

\magen{In any case, introducing flexibility in the bound to allow for arbitrary rejection sets should not be interpreted as absolving the user of any responsibility: they should still lay out the protocol they used — even if only heuristically motivated — in a convincing manner.}

\appendix

\section{Proof of Theorem~\ref{prop:unknowndep}}\label{sec:proof:prop:unknowndep}

We denote in this proof $\lambda(\alpha,X,\cH_0)$ instead of $\lambda(\alpha,\cH_0)$ to underline the dependence  of this functional w.r.t. the data $X$.
By Propositions~\ref{prop:JRadjust}
and~\ref{prop:JRadjustSD}, it sufficient to prove that $\lambda(\cdot)$ is a valid $\lambda$-calibration, that is, satisfies the requirement of Definition~\ref{def:lambdaadjust}.
Since the monotonic property is clearly satisfied, it remains to establish \eqref{lambdacalibrcontrol}.
For this, write
\begin{multline*}
  \P\bigg(\min_{1\leq k \leq K\wedge m_0}\left\{ t_k^{-1}\left(p_{(k:\cH_0)}(X)\right)\right\} < \lambda(\alpha,X,\cH_0) \bigg)\\
\begin{aligned}
&=\P\bigg( \Psi(X,\cH_0) < \lambda(\alpha,X,\cH_0) \bigg) \nonumber\\
& \leq \P\left( B^{-1}\sum_{j=1}^B \ind{ \Psi(g_j.X,\cH_0) \leq \Psi(X,\cH_0)} \leq \alpha  \right) \nonumber\\
& = \P\left( B^{-1}\sum_{j=1}^B \ind{ Y_j \leq Y_1} \leq \alpha  \right),
\end{aligned}
\end{multline*}
where we have used in the inequality the definition of $\lambda(\alpha, X,\cH_0)$ (see  \eqref{lambdaadjrand}) and we have let $Y_j=\Psi(g_j.X,\cH_0)$, $1\leq j \leq m$. Now, by \eqref{randhypo}, we easily check that $(Y_1,\dots,Y_B)$ is an exchangeable random vector: for any $g_0$ uniformly distributed on $\mathcal{G}$ (and drawn independently of the other variables),
\begin{align*}
(Y_1,\dots,Y_B)&\sim \left(\Psi(g_1.g_0.X,\cH_0),\dots,\Psi(g_B.g_0.X,\cH_0)\right)   \\
&\sim  \left(\Psi(g'_1.X,\cH_0),\dots,\Psi(g'_B.X,\cH_0)\right), 
\end{align*}
where $g'_j$, $1\leq j\leq B$, are i.i.d. uniform in $\mathcal{G}$ (independent of $X$). Above, the first equality in distribution holds because it is true conditionally on $\{g_1,\dots,g_B\}$, and the second one holds because it is true conditionally on $X$. Since the variables $\Psi(g'_j.X,\cH_0),$ $1\leq j \leq m$, are i.i.d. conditionally on $X$, we deduce that $(Y_1,\dots,Y_B)$ is an exchangeable random vector.
Hence, for any independent variable $U$ uniformly distributed on $\{1,\dots,B\}$, we obtain  
$$
\P\left( B^{-1}\sum_{j=1}^B \ind{ Y_j \leq Y_1} \leq \alpha  \right) =\P\left( B^{-1}\sum_{j=1}^B \ind{ Y_j \leq Y_U} \leq \alpha  \right).
$$
Let $\sigma$ any permutation (independent of $U$) such that $Y_{\sigma(1)}\leq \dots \leq Y_{\sigma(B)}$. Since
$\sum_{j=1}^B \ind{ Y_j \leq Y_U} = \sum_{j=1}^B \ind{ Y_{\sigma(j)} \leq Y_{U}} $
and $U$ and $\sigma(U)$ have the same distribution conditionally on $Y$, we have
\begin{align*}
&\P\left( B^{-1}\sum_{j=1}^B \ind{ Y_j \leq Y_U} \leq \alpha  \mid Y\right)
=\P\left( B^{-1}\sum_{j=1}^B \ind{ Y_{\sigma(j)} \leq Y_{\sigma(U)}} \leq \alpha \mid Y \right)\\
&\leq \P\left( B^{-1}\sum_{j=1}^B \ind{ j \leq U} \leq \alpha  \mid Y\right)
= \P\left( U \leq \alpha B  \mid Y\right) 
= \frac{\lfloor\alpha B\rfloor}{B}\leq \alpha.
\end{align*}
We underline that another argument is possible for this proof using a device
recently proposed by \cite{HG2017},
see Section~\ref{proof:improving-jr-control} for more details.

\section*{Acknowledgements}

\magen{We would like to acknowledge an associate editor and two referees for their insightful comments.}
We also thank Prof. Yoav Benjamini for interesting discussions, and Guillermo Durand for a careful reading of the manuscript.
This work has been supported by CNRS (PEPS FaSciDo), ANR-16-CE40-0019 (SansSouci) and ANR-17-CE40-0001 (BASICS).
The first author acknowledges the support from the german DFG, under the
Research Unit FOR-1735 ``Structural Inference in Statistics – Adaptation and
Efficiency'', and under the Collaborative Research Center SFB-1294 ``Data Assimilation''.

\begin{supplement} 
\sname{Supplement}
\stitle{}
\slink[doi]{COMPLETED BY THE TYPESETTER}
\sdatatype{.pdf}
\sdescription{The supplement includes relation to previous work (closed testing, higher criticism); general properties of templates and reference families; algorithms; proofs and numerical experiments. 
}
\end{supplement}

\bibliographystyle{apalike}
\bibliography{biblio}

\checknbdrafts
\end{document}